\documentclass[12pt]{amsart}
\usepackage[all,color]{xy}
\usepackage{pb-diagram}
\usepackage[mathscr]{eucal}
\usepackage{hyperref}
\hypersetup{
    colorlinks=true,
    linktoc=all,
    linkcolor=blue,
}
\usepackage{amsmath,amsfonts,amsthm,amssymb,amsthm}
\usepackage{tikz}
\usepackage{xcolor}
\usepackage{float}
\usepackage{pifont}
\usepackage{combelow}

\RequirePackage{ifthen,setspace,enumitem,tikz}
\usetikzlibrary{arrows}

\DeclareMathAlphabet{\mathpzc}{OT1}{pzc}{m}{it}

\theoremstyle{plain}

\newtheorem{lem}{Lemma}[section]
\newtheorem{cor}{Corollary}[section]
\newtheorem{theorem}{Theorem}[section]
\newtheorem{prop}{Proposition}[section]

\newtheorem{remark}{Remark}[section]

\newcommand{\pare}[1]{\left(#1\right)}
\newcommand{\brac}[1]{\left\{ #1 \right\}}
\newcommand{\sqbr}[1]{\left[ #1 \right]}
\newcommand{\norm}[1]{\left\| #1 \right\|}
\newcommand{\inner}[1]{\langle #1 \rangle}
\newcommand{\abs}[1]{\left| #1 \right|}

\newcommand{\R}{\mathbb{R}}

\newcommand{\del}{\partial}
\newcommand{\delbar}{\bar{\del}}
\newcommand{\pot}{\sqrt{-1}\del\delbar}

\renewcommand{\epsilon}{\varepsilon}

\newcommand{\vol}{\mathrm{Vol}\pare{\theta}}
\newcommand{\amp}{\mathrm{Amp}\pare{\theta}}
\newcommand{\loc}{\text{loc}}

\numberwithin{equation}{section}

\hypersetup{pdftitle={Regularity for solutions of complex Monge-Amp\`ere in big classes},
 pdfauthor={Riley Guyett, and Ryan Mc Gowan}}
\begin{document}

\title[Regularity of solutions to Monge-Amp\`ere Equations in big classes]{Regularity of solutions to complex Monge-Amp\`ere equations in big classes}
\author{Riley Guyett and Ryan Mc Gowan}
\raggedbottom
\address{Department of Mathematics, Rutgers University, Piscataway, NJ 08854}
\email{riley.guyett@rutgers.edu}
\address{Department of Mathematics, Rutgers University, Piscataway, NJ 08854}
\email{rmm396@math.rutgers.edu}

\begin{abstract}
We study the non-pluripolar complex Monge-Amp\`ere equation in an arbitrary big class on a compact K\"ahler manifold. Under a weighted Laplacian lower bound on the log density, we obtain local $C^{2,\alpha}$ regularity on the part of the ample locus where the data satisfy the stated hypotheses. A second criterion based on a quasi-plurisubharmonic decomposition gives similar regularity. In particular, a smooth positive density gives a smooth solution on the ample locus. None of these regularity results require nefness of the background class, resolving positively a question of Dinew, Guedj, and Zeriahi.
\end{abstract}

\maketitle 

\tableofcontents

\section{Introduction}

The complex Monge-Amp\`ere equation has played a central role in complex geometry since Yau's resolution of the Calabi conjecture \cite{Yau1978}, which settled the question of smooth solvability for prescribed positive volume forms in the K\"ahler class. In the case of singular data belonging to  $L^p$ for $p>1$, Bedford and Taylor developed pluripotential theoretic methods in \cite{BT} which were expanded to uniform global estimates by Ko{\l}odziej \cite{kolodziej,kolodziej2}. These fundamental results have become a central foundation of pluripotential theory and have proved to be a robust framework for the study of singular canonical metrics and degenerations that commonly arise in complex geometry; see, for instance, \cite{DNGG,DiNezzaGuedjLu,SongTian}.

The work of Boucksom, Eyssidieux, Guedj and Zeriahi \cite{BEGZ} examined the non-pluripolar Monge-Amp\`ere operator when the cohomology class $\theta$ is merely big rather than K\"ahler. There, building upon the work of the last three authors in \cite{EGZ1} examining singular K\"ahler-Einstein metrics on canonical models and klt pairs, under the additional assumption that $\theta$ is nef, they proved smoothness on the ample locus. 

Boucksom \cite{Boucksom} also characterized how the failure of nefness relates to singularities of positive currents through divisorial Zariski decompositions. Without this nefness assumption, \cite{DemaillyDinewGuedjHiepKolodziejZeriahi} established local H\"older regularity on the ample locus for smooth densities. 
Related developments include the regularity of plurisubharmonic envelopes \cite{BermanDemailly}, the regularization of closed positive currents \cite{Demailly}, and the stability of solutions in big classes \cite{GuedjZeriahi1}. A crucial missing ingredient in this regard for the regularity theory has been the establishment of complex Hessian estimates for general big classes. To this end, Dinew, Guedj and Zeriahi in \cite{DinewGuedjZeriahi} posed the following question:\medskip

\noindent \textbf{Question 21:} \emph{Let $\alpha\in H^{1,1}(X, \mathbb{R})$ be a big cohomology class on a compact manifold $X$ in the Fujiki class $\mathcal{C}$, $\mu$ a smooth volume form and let $T$ be the unique finite energy current in $\alpha$ such that $\langle T\rangle ^n=\mu$. Is $T$ smooth in $\mathrm{Amp}(\alpha)$?}\medskip

We answer in the affirmative using Corollary \ref{Answer}, combined with the reduction in Remark \ref{Fujiki}. Theorem \ref{maintheorem} treats locally H\"older log-densities on compact Kähler manifolds.

\begin{theorem}\label{maintheorem}
    Let $(X,\omega)$ be a compact K\"ahler manifold of complex dimension $n$ and let $\theta$ be a smooth closed real $(1,1)$-form representing a big class. Let $U=X\setminus Z\subseteq\amp$ with $Z$ proper analytic. Assume there exists $\phi\le -1$, smooth on $U$, such that $\theta+\pot \phi$ is a K\"ahler current having analytic singularities precisely along $Z$, and its pullback has logarithmic normal-crossing form with smooth remainder on a log resolution.
    
    Suppose $0\le f\in L^p\pare{X,\omega^n}$ for some $p>1$, $f>0$ on $U$, and $\int_X f\omega^n=\vol$. Assume $F=\log f\in C^{0,\beta}_{\loc}\pare{U}$ for some $0<\beta<1$ and satisfies either
    \begin{itemize}
        \item[(i)] $\Delta_\omega F\ge -Ce^{-B\phi}$ distributionally on $U$ for some $B,C>0$, or
        \item[(ii)] $F=F_+-F_-$ a.e. on $U$ with $F_\pm$ quasi-plurisubharmonic, $F_\pm\le C_0$, $\pot F_\pm\ge -C_0\omega$, and $F_-\in L^\infty_{\loc}\pare{U}$.
    \end{itemize}
    Then the normalized solution $\varphi\in\mathrm{PSH}\pare{X,\theta}$ of 
    \begin{equation}\label{MAEquation}
        \inner{\pare{\theta+\pot\varphi}^n} = e^F\omega^n
    \end{equation}
    belongs to $C^{2,\alpha}_{\loc}\pare{U}$ for all $0<\alpha<\beta$. If $F$ is smooth on $U$, then so is $\varphi$.
\end{theorem}

We also briefly note that the existence of a solution with minimal singularities and uniqueness up to an additive constant is provided to us by \cite{BEGZ}.

\begin{cor}\label{Answer}
    Let $\pare{X,\omega}$ be a compact K\"ahler manifold of complex dimension $n$ and let $\theta$ be a smooth closed real $(1,1)$-form representing a big class. Suppose $F\in C^\infty(X;\mathbb{R})$ satisfies $\int_X e^F\omega^n=\vol$. Then the normalized solution $\varphi\in\mathrm{PSH}\pare{X,\theta}$ to \eqref{MAEquation} is smooth on $\amp$.
\end{cor}

We explain the motivation behind conditions (i) and (ii) above a little more. Condition (i) arises naturally when pulling back a smooth positive volume form to a Kähler modification. The resulting log density has logarithmic Jacobian singularities and satisfies a uniform Laplacian lower bound. This observation allows us to extend Corollary \ref{Answer} to manifolds in Fujiki’s class $\mathcal C$; the reduction is given in Section \ref{section5}. The exponential weight in (i) permits more general deterioration of the Laplacian lower bound near $Z$.
Condition (ii) is modeled after condition $\mathcal{H}_f$ in  \cite{DinezzaLu}. There, the authors developed a capacity framework for equations like \eqref{MAEquation} with smooth data $F$ away from a Zariski closed set. Such equations arise naturally when the density may have algebraic or logarithmic singularities along the singular set; see, for instance, the works \cite{BBGZ,DDNL,DNGG,EGZ1}. 

The proof of Theorem \ref{maintheorem} consists of three main ingredients. Firstly, we consider an exhaustion of $U$ by the barrier $\phi$ in Theorem \ref{maintheorem}.  A key part of the argument is then controlling analytic estimates in terms of this background geometric description. 
Secondly, we consider a regularization of solutions of auxiliary global equations via Demailly's regularization \cite{Demailly}, coupled with the Kiselman-Legendre transform, following both the works of Eyssidieux, Guedj and Zeriahi \cite{EGZ1} and Di Nezza and Trapani \cite{DiNezzaTrapani}. Afterwards, we may apply the Richberg technique as in \cite{HLP} to regularize subsolutions. 
Finally, motivated by the classical work of \cite{CKNS}, we combine the geometric estimates with our regularized subsolutions to reduce the problem to solving a Dirichlet problem via the work of Guan and Li \cite{GuanLi}. Here, the geometric and regularization estimates come together: the smooth strict subsolutions supply boundary data for the Dirichlet problems, while the control of the exhaustion geometry allows quantitative estimates for their solutions. With pointwise gradient estimates following \cite{PS,PS2} and trace estimates in place, we may then pass to the limit and establish the claimed regularity.

The organization of this article is as follows. Section \ref{section2} regularizes the densities and constructs an exhaustion of $U$ on which we can establish H\"older estimates for solutions of auxiliary equations. Section \ref{section3} produces smooth strict subsolutions with controlled approximation error and fourth-derivative estimates. Section \ref{section4} solves a family of Dirichlet problems with these subsolutions as boundary data, and provides weighted gradient estimates along the exhaustion. Section \ref{section5} extends the gradient estimates to boundary Hessian estimates and interior trace estimates. We thus obtain local uniform ellipticity. Identifying the limiting behavior of our Dirichlet solutions and our global weak solutions, we verify that their common limit solves the Monge-Amp\`ere equation with the desired regularity. 

\bigskip
\noindent\textbf{Acknowledgments.}
The authors thank their advisor, Professor Jian Song, for introducing them to this problem and for his encouragement. This work forms part of the authors' Ph.D. theses at Rutgers University, and they are grateful to the Department of Mathematics for its generous support. \bigskip 

\noindent \textbf{AI Declaration:} The authors’ starting point was an experimental idea: to approximate the global equation by solutions of Dirichlet problems for the complex Monge-Ampère equation, using a quantitative analytic exhaustion of the Kähler locus. The proof was ultimately completed by ChatGPT Sol on September 8, 2026, in the case when the right-hand side volume measure is positive and smooth. The authors have since revised the argument and independently verified all results, and they take full responsibility for the content of this paper.

\section{Preliminaries and Initial Estimates}\label{section2}

Let $(X,\omega)$ be a compact K\"ahler manifold of complex dimension $n$. Let $\theta$ be a smooth closed $(1,1)$-form representing a big class. Define the ample locus $\amp$ to be the collection of all points near which some K\"ahler current in $\sqbr{\theta}$ with analytic singularities is smooth.

We say that a function $\varphi:X\to \mathbb{R}\cup\{-\infty\}$ is quasi-plurisubharmonic if locally we can write $\varphi=\varphi_1+\varphi_2$ where $\varphi_1$ is smooth and $\varphi_2$ is plurisubharmonic. We furthermore say that $\varphi$ is $\theta$-plurisubharmonic if $\varphi$ is quasi-plurisubharmonic and
\[\theta_\varphi := \theta+\pot \varphi\geq 0\] 
in the weak sense of currents. We denote the space of all $\theta$-plurisubharmonic functions on $X$ by $\mathrm{PSH}(X,\theta)$. We say $\varphi\in \mathrm{PSH}(X,\theta)$ has minimal singularities if $\varphi-V_\theta$ is bounded, where $V_\theta$ is the envelope
\[
V_\theta=(\sup\{\varphi\in \operatorname{PSH}(X,\theta): \varphi\leq 0\})^*.
\]

If $Z=\varnothing$ and $F$ is smooth, then $\theta$ lies in a K\"ahler class and the regularity of the solution follows from Yau's theorem \cite{Yau1978}.

\subsection{Barrier Estimates along an Exhaustion} Let $\phi$ be as in the statement of Theorem \ref{maintheorem}. This implies the existence of $\delta_0$ such that 
\begin{equation}
    \theta_\phi \ge 2\delta_0\omega.
\end{equation}
Choose $\sigma>0$ sufficiently small so that $\psi=\pare{1+\sigma}\phi$ satisfies $\theta_\psi\ge \delta\omega$ for some $\delta>0$. Therefore, exponential bounds in $-\phi$ can be taken to be exponential bounds in $-\psi$. We will later subtract a fixed constant from $\phi$ and adjust $\psi$ and $\sigma$ accordingly so that $\psi=\pare{1+\sigma}\phi$ still holds; see Remark \ref{barriershift}. This shift will not affect the following estimates.

For each positive integer $j$, we define the set 
\begin{equation}
    \Omega_j = \brac{\psi>-j}\Subset U.
\end{equation}
Since $\psi$ also has singularities precisely along $Z$ and tends to $-\infty$ along $Z$, the sets $\brac{\Omega_j}$ satisfy $\overline{\Omega}_j\Subset\Omega_{j+1}$ and form an exhaustion of $U$. The next few lemmata describe the interplay between the regularity of the barrier $\psi$ and the geometry of the sets $\Omega_j$.

\begin{lem}\label{psibounded}
For each positive integer $j$, on $\Omega_{j+1}$, we have 
\begin{equation*}
    \abs{\nabla^k\psi}_\omega \le C_ke^{B_kj} \quad \pare{k\ge 1}.
\end{equation*}
Moreover, on $U$ we have a pointwise bound
\begin{equation*}
    \abs{\nabla^k\psi}_\omega \le C_ke^{-B_k\psi} \quad \pare{k\ge 1}.
\end{equation*}
\end{lem}

\begin{proof}
We resolve the singularities of $\phi$ as in \cite{DiNezzaTrapani}. Fix a log-resolution $\pi\colon \hat{X}\to X$, biholomorphic over $U$. Choose a finite collection of relatively compact normal-crossing charts $\{U_\nu\}$ with coordinates $z=(z_1, \dots, z_n)$, on each of which we have
\begin{equation}
    \pi^*\phi(z) = \sum_{\ell=1}^r a_\ell^{(\nu)}\log\abs{z_\ell}^2 + h_\nu(z)
\end{equation}
where $D\cap U_\nu=\{z\in \hat{X}: z_1\cdots z_{r_\nu}=0\}$, $a_\ell^{(\nu)}>0$ for $\ell=1, \dots, r_\nu$, and $h_\nu(z)\in C^\infty(U_\nu)$.
In particular, each chart $U_\nu$ is relatively compact in a larger normal-crossing chart $U_\nu'$ to which $h_\nu(z)$ extends. Thus, for each fixed integer $k\geq 0$ we have
\[
|D^kh_\nu(z)| \leq H_k
\]
uniformly on each smaller chart $U_\nu$. 

Fix a smooth metric on the compact resolution $\hat{X}$. In each fixed chart on $\hat{X}$, suppressing $\nu$ for simplicity, we have
\begin{equation*}
    \pi^*\psi(z) = \sum_{\ell=1}^r (1+\sigma)a_\ell\log\abs{z_\ell}^2 + (1+\sigma)h(z).
\end{equation*}
After rescaling the fixed charts $U_\nu$ to some $\hat{U}_\nu$, we may assume that $\abs{z_\ell}<1$ for all $\ell=1,\dots, n$. The remainder $(1+\sigma)h(z)$ and each of its fixed-order derivatives are uniformly bounded on each $\hat{U}_\nu$. If $\psi>-j-1$, then 
\begin{equation*}
    \sum_{\ell=1}^r (1+\sigma)a_\ell\log\abs{z_\ell}^2 \ge -j-C.
\end{equation*}
This forces every non-zero divisor coordinate to be bounded below by $ce^{-Bj}$ on $\pi^{-1}\pare{\Omega_{j+1}}$. Then for each $k\ge 1$, differentiation yields
\begin{equation*}
    \abs{D^k\pi^*\psi}_{\hat{\omega}}\le C_ke^{B_kj},
\end{equation*}
on $\pi^{-1}(\Omega_{j+1})$, where $\hat{\omega}$ is a fixed smooth hermitian metric on $\hat{X}$.

To transfer these estimates to $X$, we use coordinate representations of $\pi$ between the fixed charts. It is biholomorphic away from the divisor, and its holomorphic Jacobian (which vanishes only along that divisor) can be written locally as
\begin{equation*}
    \det(d\pi(z)) = u(z)\prod_{\ell=1}^rz_\ell^{q_\ell},
\end{equation*}
where $u(z)$ is a nonvanishing holomorphic function on the larger chart $\hat{U}_\nu$, and $q_\ell\ge 0$. On the smaller chart, $u^{-1}(z)$ and its derivatives are bounded. The adjugate formula gives
\[
(d\pi)^{-1}(z)=\frac{\mathrm{adj}(d\pi(z))}{\det(d\pi(z))}.
\]
In particular, the coefficients of $\pare{d\pi}^{-1}(z)$ and their fixed-order derivatives have at most finite-order poles along the divisor. The chain rule then provides
\begin{equation*}
    \abs{\nabla^k\psi}_{\omega}\le C_ke^{B_kj}\qquad \text{on }\Omega_{j+1},\quad k\ge 1.
\end{equation*}
For arbitrary $x\in U$, we can then consider the smallest integer $j_0\geq 1$ for which $x\in \Omega_{j_0+1}$. If $j_0>1$, by minimality we have that $x\notin \Omega_{j_0}$ so $\psi(x)\leq -j_0$. In particular, we find $j_0\leq -\psi(x) +C$, including the case where $j_0=1$. Hence
\[
|\nabla^k\psi(x)|_\omega\leq C_ke^{B_k j_0}\leq C_ke^{-B_k \psi(x)}.
\qedhere\]
\end{proof}

\begin{lem}[]\label{exhaustioncovering}
For sufficiently large integers $j$, there are fixed constants $c,C>0$ such that $\overline{\Omega}_j$ admits a cover by at most $N_j:=Ce^{Bj}$ holomorphic coordinate balls $B^{(\nu)}_{ce^{-bj}}\pare{x_\nu}$ for $\nu=1,\dots, N_j$ such that $B^{(\nu)}_{3ce^{-bj}}(x_j)\Subset \Omega_{j+1}$.
\end{lem}

\begin{proof}
    A path $\gamma$ from $\overline{\Omega}_j$ to $X\setminus \Omega_{j+1}$ must cross from the level set $\psi=-j$ to the level set $\psi=-j-1$. On the band $\Omega_{j+1}\setminus\overline{\Omega_j}$
    by Lemma \ref{psibounded} we have 
        \begin{equation*}
            \abs{\nabla\psi}_{\omega} \le Ce^{Bj}.
        \end{equation*}
        Integrating along the crossing portion of the path, we see
        \begin{equation*}
            Ce^{Bj}\operatorname{Length}_\omega(\gamma)\geq 1.
        \end{equation*}
        Taking the infimum, we have 
        \begin{equation}
            \operatorname{dist}_\omega(\Omega_j,X\setminus \Omega_{j+1}) \ge ce^{-bj}.
        \end{equation}
        Taking $R_j=ce^{-bj}$, after perhaps shrinking $c$, we see that every coordinate ball of radius $3R_j$ centered in $\Omega_j$ is compactly contained in $\Omega_{j+1}$. A maximal separated set of centers in $\overline{\Omega}_j$ gives a covering. Because $X$ is of real dimension $2n$, the number of balls needed is bounded by 
        \begin{equation*}
            CR_j^{-2n} \le N_j.
            \qedhere
        \end{equation*}
\end{proof}

\begin{remark}
    The purpose of the exhaustion sets $\{\Omega_j\}$ is twofold. Firstly, as seen above, the sets $\Omega_j$ provide a convenient covering of $U$. As a result, the arguments in Section \ref{section3} can be localized to holomorphic coordinate balls. Secondly, we shall see in Section \ref{section4} that the geometry of these sets is well-behaved vis-\`a-vis solving certain Dirichlet problems. We defer proving these geometric properties to Section \ref{section4}.
\end{remark}

\subsection{Regularization of the Density}

In Theorem \ref{maintheorem}, we are merely assuming that $f\in L^p(X,\omega^n)$, and as a result, it may possess singular behavior. We may regularize the density $f$ by utilizing heat kernels and heat operators. Since these heat kernel estimates are only needed for this subsection, we refer the interested reader to Appendix \ref{heatappendix} for the details.  

\begin{lem}\label{Fapprox}
    There are smooth positive densities $f_j=e^{F_j}$ on $X$ with mass $\vol$, with the following properties:
    \begin{enumerate}
        \item The $f_j$ are uniformly bounded in $L^p\pare{X,\omega^n}$ and $f_j\to f$ in $L^1\pare{X,\omega^n}$. The $F_j$ are uniformly bounded in $C^{0,\beta}_{\loc}\pare{U}$ and $F_j\to F$ locally uniformly on $U$.
        \item On $X$, $ce^{-j}\le f_j\le Ce^{Bj}$ and $\abs{\nabla F_j}_\omega\le Ce^{Bj}$.
    \end{enumerate}
    Further, under the respective hypotheses of Theorem \ref{maintheorem}:
    \begin{itemize}
        \item[(i)] On $\Omega_{j+1}$, $F_j\le C-B_1\psi$ and $\Delta_\omega F_j\ge -Ce^{-B_2\psi}$.
        \item[(ii)] There is a smooth decomposition $F_j=F_j^+-F_j^-$ with $F_j^\pm\le C_*$, $\pot F_j^\pm\ge -C_*\omega$. Moreover, for every $K\Subset U$, there is a constant $C_K>0$ independent of $j$ such that $F_j^-\ge -C_K$ on $K$.
    \end{itemize}
\end{lem}

Although both cases use heat regularization, different constructions are needed to obtain properties (i) and (ii). We therefore split the proof into two parts.

\begin{proof}[Proof under hypothesis (i)]
    Fix $x\in U$, and choose $k$ with $k=-\psi\pare{x}+O(1)$ and $x\in \Omega_k$. Choose $r_k=ce^{-bk}$ so that $B_{2r_k}\pare{x}\Subset\Omega_{k+1}$. On this ball, $\Delta_\omega F\ge -Ce^{B\pare{k+1}}$ distributionally. Increasing $b$ so that $r_k^2Ce^{B\pare{k+1}}\le 1$ and Lemma \ref{chainrule} imply $\Delta_\omega e^F\ge -Ce^{B\pare{k+1}}e^F$ in distribution.
    The local mean-value inequality for nonnegative subsolutions yields
    \[
    e^{F\pare{x}}\le Cr_k^{-2n}\int_{B_{r_k}\pare{x}} e^F\omega^n\le Ce^{Bk}.
    \]
    Since $k=-\psi\pare{x}+O\pare{1}$, we see that $e^F\le Ce^{-B\psi}$ on $U$.

    On $U$, set
    \[
    h_j\pare{s} = \log \pare{e^s+e^{-j}},\qquad G_j=\min\brac{h_j\pare{F},M_j}
    \]
    where $M_j=C_1\pare{j+3}$. We see that $G_j$ satisfies
    \[
    -j\le G_j\le M_j,\qquad e^{G_j}\le e^F+e^{-j}
    \]
    almost everywhere. Extending $G_j$ to be $-j$ on $X\setminus U$ extends these bounds to hold throughout $X$. Our choice of $M_j$ makes the cap inactive on $\Omega_{j+3}$, so on that exhaustion level
    \[
    G_j=h_j\pare{F},\qquad \Delta_\omega G_j\ge -Ce^{-B\psi},
    \]
    with the Laplacian bound following from Lemma \ref{chainrule}.

    Setting $t_j=t_0e^{-Dj}$ for $D>2b$ with $b$ as above, we define the regularization
    \[
    \tilde{G}_j=P_{t_j}G_j,\qquad c_j=\frac{\vol}{\int_X e^{\tilde{G}_j}\omega^n},\qquad F_j=\tilde{G}_j+\log c_j
    \]
    where $P_t$ is the scalar heat operator. 

    Fix $K\Subset V\Subset U$ with $V$ open and choose $\eta\in C^\infty_c\pare{U}$ equal to one on a neighborhood of $\overline{V}$. Set $g=\eta F$ and $g_j=\eta h_j\pare{F}$, extended by zero to be defined on all of $X$. On $\operatorname{supp}\eta$,
    \[
    0\le h_j\pare{F}-F = \log\pare{1+e^{-j-F}}\le e^{-j-\underset{\mathrm{supp}\eta }{\inf}F}.
    \]
    Therefore, $g_j\to g$ uniformly and because $0\le h_j^\prime\le 1$, the $g_j$ are uniformly bounded in $C^{0,\beta}\pare{X}$. For $j$ sufficiently large $G_j=g_j$ on $V$. By our choice of $D$ and Lemma \ref{heatconverge} we have $\tilde{G}_j\to F$ locally uniformly and each $\tilde{G}_j$ has uniform $C^{0,\beta}_{\loc}$ bounds. Since the analytic set $Z$ has measure $0$, we have that $e^{\tilde{G}_j}\to e^F$ almost everywhere. Jensen's inequality for the probability kernel of $P_t$ followed by the invariance of $\omega^n$ implies that each $e^{\tilde{G}_j}$ is uniformly bounded in $L^p$. By Vitali's convergence theorem, we conclude that $e^{\tilde{G}_j}\to e^F$ in $L^1$. In particular, $c_j\to 1$, proving item (1).

    To prove item (2), Lemma \ref{heatbounded} yields
    \[
    \norm{\nabla F_j}_{L^\infty\pare{X,\omega^n}} \le Ct_j^{-1/2}\norm{G_j}_{L^\infty\pare{X,\omega^n}}\le Ce^{Bj}.
    \]
    Thus, $-j\le \tilde{G}_j\le M_j$.

    Fix $x\in \overline{\Omega}_{j+1}$ and let $k$ be the smallest integer strictly greater than $-\psi\pare{x}$. Then $x\in \Omega_k$, $k\le j+2$, and $k=-\psi+O\pare{1}$. 
    Consider a smooth cutoff function $\chi\in C^\infty_c(B_{2r_k}(x))$ satisfying $\mathrm{supp}(\chi)\subseteq B_{2r_k}(x)$, $0\leq \chi\leq 1$, $\chi\equiv 1$ on $B_{r_k}(x)$ and such that $|D\chi|\leq Cr_{k}^{-1}$ and $|D^2\chi|\leq Cr_k^{-2}$.
    Decomposing $\Delta_\omega F_j\pare{x}$ as 
    \[
    \Delta_\omega F_j\pare{x} = \inner{\Delta_\omega G_j,\chi K_{t_j}\pare{x,\cdot}}+\int_X G_j\pare{y}\Delta_{\omega,y}\sqbr{\pare{1-\chi\pare{y}}K_{t_j}\pare{x,y}}\omega^n\pare{y},
    \]
    we see that each term of $\Delta_{\omega,y}\big((1-\chi(y))K_{t_j}(x,y)\big)$
    is supported at a distance of at least $r_k$. If we call this integral $R_j\pare{x}$, Lemma \ref{heatbounded} gives the following bound:
    \[
    \abs{R_j\pare{x}}\le C\pare{j+1}r_k^{-2}t_j^{-N}e^{-cr_k^2/t_j}.
    \]
    Since $D>2b$, we see that $R_j(x)$ is uniformly bounded and furthermore that $R_j(x)\to 0$ as $j\to \infty$. On $\mathrm{supp}(\chi)\subset\Omega_{k+1}\subset \Omega_{j+3}$, $\chi K_{t_j}\ge 0$ and $\inner{\Delta_\omega G_j,\chi K_{t_j}\pare{x,\cdot}}\ge -Ce^{B\pare{k+1}}$. Hence, for $k=-\psi\pare{x}+O(1)$, we have
    \[
    \Delta_\omega F_j(x)\geq -Ce^{B(k+1)}.
    \]

    On $B_{r_k}\pare{x}$, our preliminary growth estimate gives $G_j\le C\pare{1+k}$. Outside the ball, $G_j\le M_j$, and the heat kernel mass is at most $Ct_j^{-n}e^{-cr_k^2/t_j}$. Thus,
    \[
    \tilde{G}_j\le C\pare{1+k}+CM_jt_j^{-n}e^{-cr_k^2/t_j} \le C-B_1\psi\pare{x}.
    \]
    Adding the bounded $\log c_j$ completes claim (i) and the proof.
\end{proof}

\begin{proof}[Proof under hypothesis (ii)]
    Since $Z$ is analytic and the summands are bounded above, extend $F_\pm$ across $Z$ by the removable singularity theorem as quasi-psh functions preserving the complex Hessian lower bound. If $G_\pm$ are the extensions, then $G_\pm\in L^1\pare{X,\omega^n}$ and $F=G_+-G_-$ a.e. We also have $\pot G_\pm\ge -C_0\omega$, so if $P_t$ denotes the scalar heat operator, then $G_t^\pm=P_tG_\pm$ satisfies $\pot G_t^\pm\ge -C_1\omega$ by Lemma \ref{heathessian}.

    If $t_j=t_0e^{-j}$, we write $G_j^{\pm}$ for $G_{t_j}^\pm$ and define $H_j=G_j^+-G_j^-$. We claim the approximations are
    \begin{equation}
        f_j = c_j\pare{e^{H_j}+e^{-j}},\qquad F_j=\log f_j,
    \end{equation}
    where $c_j$ is a normalization constant to preserve mass. Fix $K\Subset V\Subset U$ with $V$ open and choose $\eta\in C^\infty_c\pare{U}$ equal to one on a neighborhood of $\overline{V}$. Taking $a_j=G_+-G_-$ and $g_j=g=\eta F$, Lemma \ref{heatconverge} shows that $H_j$ converges locally uniformly to $F$ and that the $H_j$ have uniformly bounded local $C^{0,\beta}$ norms. Jensen's inequality for the measure $K_{t_j}\pare{x,y}\, \omega^n\pare{y}$ gives a uniform $L^p$ bound for $e^{H_j}$, and Vitali gives $e^{H_j}+e^{-j}\to f$ in $L^1$, so $c_j\to 1$, hence proving item (1).
    
    The kernel bounds of Lemma \ref{heatbounded} yield
    \[
    e^{H_j}\le C t_j^{-n/p}\norm{f}_{L^p},\quad \abs{\nabla F_j}_\omega = \frac{e^{H_j}}{f_j/c_j}\abs{\nabla H_j}_\omega \le Ct_j^{-n-1/2}\norm{G_+-G_-}_{L^1},
    \]
    which proves item (2).

    To see that $F_j$ decomposes, let 
    \begin{equation}
        F_j^-=G_j^-,\qquad F_j^+=\log\pare{e^{G_j^+}+e^{G_j^- -j}}+\log c_j.
    \end{equation}
    We see that indeed $F_j=F_j^+-F_j^-$ and both summands are bounded above uniformly by the scalar maximum principle and the fact that $c_j\to 1$. The fact that $\pot F_j^-\ge -C_*\omega$ follows from Lemma \ref{heathessian}.
    Because the complex Hessian of $\log\pare{e^a+e^b}$ is a convex combination of the Hessians of $a,b$, and a semipositive term, $\pot F_j^+\ge -C_*\omega$ as well.

    For $K\Subset V\Subset U$, with $V$ chosen so that $G_-\ge -C_V$ a.e., we split the heat integral defining $F_j^-\pare{x}$ for $x\in K$ into the contribution from $V$ and the contribution from $X\setminus V$. The contribution from $V$ is at least $-C_V$ since the scalar kernel is positive and of mass one. The contribution from $X\setminus V$ has absolute value bounded by $Ct_j^{-n}e^{-c_K/t_j}\norm{G_-}_{L^1}$. Hence, $F_j^-\ge -C_K$ uniformly on $K$, proving item (ii) and completing the proof.
\end{proof}

\begin{cor}\label{smoothFapprox}
    In the setting of Lemma \ref{Fapprox}, if $F$ is smooth, then $F_j\to F$ in $C^\infty_{\loc}\pare{U}$.
\end{cor}
\begin{proof}
    Fix $K\Subset V\Subset U$ with $V$ open and choose $\eta\in C^\infty_c\pare{U}$ equal to one on a neighborhood of $\overline{V}$. Extend all products with $\eta$ to be zero in $X$. We set
    \[
    \pare{a_j,g_j} = \begin{cases}
        \pare{G_j,\eta h_j\pare{F}} & \text{in case (i)} \\
        \pare{G_+-G_-,\eta F} & \text{in case (ii)}.
    \end{cases}
    \]
    We see that $a_j=g_j$ on $V$ almost everywhere for sufficiently large $j$. 

    We claim that $P_{t_j}a_j\to F$ in $C^\infty\pare{K}$. To see this, decompose
    \[
    P_{t_j}a_j = P_{t_j}g_j+P_{t_j}\pare{a_j-g_j}.
    \]
    For the first term, $P_{t_j}g_j\to \eta F$ in $C^\infty\pare{X}$. In case (ii) this is immediate because $g_j$ is constant in $j$. Case (i) is more subtle. Since $F$ is smooth on $\operatorname{supp}\eta$, its image is contained in some compact interval $I\subset\R$. On $I$, $h_j\pare{s}\to s$ smoothly, since 
    \[
    h_j\pare{s} = s+\log\pare{1+e^{-j-s}}
    \]
    and every derivative on the right-hand side tends uniformly to zero on $I$. Therefore, $g_j\to \eta F$ in $C^\infty\pare{X}$. The scalar heat operators are bounded for small time in each fixed $C^m$ norm by Lemma \ref{heatbounded} and tend to the identity in these norms on smooth functions. Therefore,
    \[
    P_{t_j}g_j-\eta F = P_{t_j}\pare{g_j-\eta F}+\pare{P_{t_j}\eta F-\eta F}\to 0\quad\text{in }C^m\pare{X}
    \]
    for every $m\ge 0$.

    For the second term, put $\rho=\operatorname{dist}\pare{K,X\setminus V}>0$. The off-diagonal heat kernel estimates up to order $m$ give
    \[
    \norm{P_{t_j}\pare{a_j-g_j}}_{C^m\pare{K}}\le C_mt_j^{-n-m/2}e^{-c_m\rho^2/t_j}\norm{a_j-g_j}_{L^1\pare{X}}.
    \]
    The last factor is $O\pare{j+1}$ in case (i) and bounded in case (ii). Since $t_j$ decreases exponentially, the right-hand side tends to zero. Since $\eta F=F$ on $K$, we are done.
\end{proof}

Thanks to the work of Boucksom, Eyssidieux, Guedj and Zeriahi \cite[Theorem 3.1, Theorem 4.1]{BEGZ}, there exists a unique solution $\varphi$ with minimal singularities to the following complex Monge-Amp\`ere equation
\begin{equation*}
\langle(\theta+i\partial\overline{\partial} \varphi)^n\rangle =e^{F}\omega^n.
\end{equation*}
We now develop a family of auxiliary equations whose solutions will converge to $\varphi$. If $\lambda_j=e^{-j}$, then by \cite[Theorem 4.1, Theorem 6.1]{BEGZ}, we can solve the corresponding equation in the class $\lambda_j\sqbr{\theta}$, with smooth positive volume form $\lambda_j^nf_j\omega^n$, and divide by $\lambda_j$ to find solutions $\varphi_j$ to
\begin{equation}\label{BEGZj}
\langle(\theta+i\partial\overline{\partial} \varphi_j)^n\rangle =e^{\lambda_j\varphi_j+F_j}\omega^n.
\end{equation}

\subsection{Regularization of the Solutions to the Auxiliary Equations}

We now turn our attention to deriving estimates on the potentials $\varphi_j$ uniformly in $j$.

\begin{lem}\label{unifpotentialboundi}
    There is a uniform constant $C>0$ such that for all $j$,
    \[
    V_\theta-C\leq \varphi_j\leq C.
    \]
\end{lem}

\begin{proof}
    For the upper bound, let $m_j =\sup_X \varphi_j$, and consider 
    \begin{equation}\label{difference}
    \psi_j = \varphi_j -m_j.
    \end{equation}
    Integrating \eqref{BEGZj} we find \[
    \vol = \int_X \inner{(\theta+i\partial\overline{\partial} \varphi_j)^n} = e^{\lambda_j m_j}\int_X e^{\lambda_j \psi_j+F_j}\omega^n \le e^{\lambda_j m_j}\vol.
    \]
    Since $\psi_j$ is a normalized $\theta$-psh function, by the uniform exponential integrability of \cite[Theorem 4.1]{BEGZ}, we see $\psi_j $ has a uniform $L^{p/\pare{p-1}}$-bound. Furthermore, $\psi_j\leq 0$. Hence, applying Jensen's inequality, we find
    \begin{align*}
    0 &\le m_j\le \frac{1}{\vol}\int_X-\psi_j e^{F_j}\,\omega^n\le C\norm{\psi_j}_{L^{p/\pare{p-1}}}\norm{e^{F_j}}_{L^p}\le C.
    \end{align*}
    For the lower bound, the normalized potential $\psi_j$ has density $e^{\lambda_j\varphi_j}f_j$, which is uniformly bounded in $L^p$ by the upper bound just proved for $\varphi_j$ and Lemma \ref{Fapprox}. The relative estimate of \cite[Theorem 4.1]{BEGZ} and $m_j\ge 0$ hence provide
    \[
    V_\theta-C\leq \psi_j\le \varphi_j\le C.\qedhere
    \]
\end{proof}

\begin{remark}\label{barriershift}
   By Lemma \ref{unifpotentialboundi}, we can subtract a fixed constant from $\phi$ and $\psi$ so that $\phi,\psi\leq \psi_j$ for all $j$ while retaining $\psi=\pare{1+\sigma}\phi$.  
\end{remark}

The goal now is to regularize $\psi_j$, which in turn regularizes $\varphi_j$. Consider Demailly's regularization operator $\Psi(u)$ defined by \cite{Demailly,DiNezzaTrapani}
\[
\Psi(u)(z,w) = \int_{\zeta \in T_zX} u(\operatorname{exph}_z(w\zeta))\chi(|\zeta|^2)\, dV_\omega \quad \text{for }(z,w)\in X\times \mathbb{C}
\]
where $\operatorname{exph}:TX\to X$ is the formal holomorphic part of the Taylor series expansion of the exponential map of the Chern connection on $TX$ associated to the metric $\omega$, $\abs{w}$ is taken sufficiently small so that the exponential map is defined, and $\chi:\mathbb{R}\to \mathbb{R}^+$ is a smooth cut-off function such that $\chi>0$ for $t<1$, $\chi(t)=0$ for $t\geq 1$ and $\int_{\mathbb{C}^n} \chi(|z|^2)\, d\lambda(z)=1$ where $d\lambda(z)$ denotes the usual Lebesgue measure on $\mathbb{C}^n$.

Note that if $w=0$, one has $\Psi(u)(z,0)=u(z)$. Furthermore, for $\delta>0$, there is a constant $K>0$ (depending on curvature) such that for $|w|\leq \delta$ we have 
\begin{equation}\label{curvest}
\Psi(u)(z,w)-u(z)\geq -K\abs{w}^2.
\end{equation}
Setting $\abs{w}=\delta$, we define the Kiselman-Legendre transform on $X$ with level $c>0$ to be 
\[
U_{c,\delta}(z)=\inf_{0<t\leq\delta} \left\{\Psi(\psi_j)(z,t) +K(t^2-\delta^2)+K(t-\delta)-c\log\left(\frac{t}{\delta}\right)\right\}.
\]
By the work of \cite[Equation 1.13]{BermanDemailly} and \cite[Lemma 4.1]{KolodziejNguyen}, there exists some constant $A$ such that
\begin{equation}\label{diNezzaTrapaniA}
\theta+i\partial\overline\partial U_{c,\delta}\geq -(Ac+2K\delta)\omega
\end{equation}
and furthermore 
\begin{equation}\label{diNezzaTrapaniB}
\psi_j(z)-K(\delta+\delta^2)\leq U_{c,\delta}(z)\leq \Psi(\psi_j)(z,\delta)\leq 0.    
\end{equation}
The constants $A,K$ are seen to be independent of $j$ by using the Jensen-Lelong monotonicity formula. 

\begin{lem}[]\label{Psi-psiL1bound}
    Let $\psi_j$ be defined as in \eqref{difference}. The regularization obeys the following bound for $|w|=\delta$
    \[
    \|\Psi(\psi_j)(z,w)-\psi_j(z)\|_{L^1(X,\omega)} \leq C\delta,
    \]
    independent of $j$.
\end{lem}

\begin{proof}
    The idea is to prove the estimate locally on holomorphic coordinate balls $B'\Subset B$ which form an atlas for $X$, and then patch the estimates together via compactness of $X$. 
    To begin, as $\psi_j \in \mathrm{PSH}(X, \theta)$, we have 
    \[
    \Delta_\omega\psi_j \geq -C.
    \]
    Also, since $\omega$ is K\"ahler, we have
    \[
    \int_X (\Delta_\omega \psi_j) \omega ^n = n\int_X (i\partial \overline\partial \psi_j)\wedge \omega^{n-1}=0.
    \]
    Combining these facts, we conclude that the measure
    \[
    \nu_j = (\Delta_\omega \psi_j)dV_\omega
    \]
    has bounded total variation. 
    Secondly, the Green's representation formula gives the $L^1$-bound 
    \begin{equation}\label{l1bound}
    \|\nabla\psi_j\|_{L^1(X)}\leq C.
    \end{equation}
    We now focus on handling the regularization. If we consider paths 
    \[
    \gamma:[0,1]\to X, \qquad \gamma(\tau) = \mathrm{exph}_z(\tau w\zeta)
    \]
    then $|\gamma'(\tau)|\leq C'|w|$. The following fundamental theorem of calculus argument is first carried out for smooth functions. By the density of $C^\infty(X)$ in $W^{1,1}(X)$, and the $L^1$-boundedness of the regularization operator, the estimate extends to $\psi_j$
    \[
    \int_{B'}|\psi_j (\mathrm{exph}_z(w\zeta))-\psi_j (z)|\, dV_\omega(z)
    \leq C'|w|\int_0^1 \int_{B'}|\nabla \psi_j(\gamma(\tau))|\, dV_\omega(z)\, d\tau.
    \]
    By \eqref{l1bound}, we may apply Fubini's theorem to swap the order of integration above. By enlarging to $B\Supset B'$, and performing a change of coordinates, we get  
    \begin{align*}
    \int_{B'}|\psi_j (\mathrm{exph}_z(w\zeta))-\psi_j (z)|\, dV_\omega(z)
    &\leq C'|w|\int_B |\nabla\psi_j (z')|\, dV_\omega(z')\\
    &\leq C'\delta \|\nabla \psi_j\|_{L^1(B)}\leq C\delta.
    \end{align*}
    Integrating against the cut-off function and using the preceding approximation argument gives the estimate
    \begin{align*}
    \|\Psi(\psi_j)(z,w)-&\psi_j(z)\|_{L^1(X,\omega)}\\
    &\leq \int_X\int_{|\zeta|\leq 1} |\psi_j(\mathrm{exph}_z(w\zeta))-\psi_j(z)|\chi(|\zeta|^2)\, dV_\omega(\zeta)\, dV_\omega(z).\qedhere
    \end{align*}
\end{proof}

We need the following stability estimate, whose formulation comes from \cite[Proposition 2.9]{DemaillyDinewGuedjHiepKolodziejZeriahi} but was originally proved in \cite{GuedjZeriahi1}.
\begin{prop}\label{GZProp}
    Assume that $\mu$ is a probability measure absolutely continuous with respect to a smooth volume form $dV$, $d\mu = f dV$, where $f \in L^p(X)$ with $p > 1$. Let $\varphi,\psi\in\mathrm{PSH}\pare{X,\theta}$ satisfy
    \[
    \mathrm{MA}(\varphi) = \mu, \qquad  -M_0 + V_\theta \leq  \varphi \leq V_\theta, \qquad  \psi \leq  V_\theta\]
    for some positive constant $M_0>0$. Then for any exponent $0 < \gamma < \frac{p-1}{np+p-1}$, there exists a constant
    $B_0 = B_0(p, \gamma, M_0,\norm{f}_{L^p},\theta) > 0$ such that
    \[
    \sup_X(\psi-\varphi)_+ \leq  B_0\|(\psi-\varphi)_+\|^\gamma_{
    L^1(X)}.
    \]
\end{prop}
\begin{remark}
    We note that $\mu$ being a probability measure is simply a normalization convention. By rescaling the reference form $\theta$ and potentials, the same estimate holds for $\mu(X)=\mathrm{Vol}(\theta)$. All that changes is that the quantity $B_0$ may additionally depend on the total measure $\mu(X)$.
\end{remark}

We turn the $L^1$ estimate for Demailly regularization into pointwise control.

\begin{lem}\label{Holdercont}
    Let $\psi_j$ be defined as in \eqref{difference}. For $r=\abs{w}$ sufficiently small, there exists $0<\beta_0<1$ such that 
    \[
    |\Psi(\psi_j)(z,w)-\psi_j(z)|\leq Ce^{BH(z)}r^{\beta_0},
    \]
    where $H(z)=V_\theta(z)-\phi(z)$.
\end{lem}

\begin{proof}
    This proof strategy is based on \cite{DemaillyDinewGuedjHiepKolodziejZeriahi}, with the linear correction term as noted in \cite{DiNezzaTrapani}. Fix $\gamma$ as in Proposition \ref{GZProp} and let $0<\beta_0<\min\brac{\gamma,1}$. For $0<s<s_0$ small, let $\eta=s^{\beta_0}$ and $c_s=\delta_0\eta/A$. Define the convex combination
    \begin{equation}\label{convexcombo}
    W_s=\pare{1-\eta}U_{c_s,s}+\eta\phi,
    \end{equation}
    where $s_0\ll1$ is sufficiently small so that $2Ks\le \delta_0\eta$. By \eqref{diNezzaTrapaniA}, we have
    \begin{align*}
        \theta_{W_s} &\ge \pare{2\delta_0\eta-\pare{1-\eta}\pare{Ac_s+2Ks}}\omega\ge 0.
    \end{align*}
    Since $W_s\leq 0$, we conclude that $W_s\le V_\theta$. Since $\phi\le \psi_j$, by Lemma \ref{Psi-psiL1bound},
    \begin{align*}
        \norm{\pare{W_s-\psi_j}_+}_{L^1}\le \norm{\pare{\Psi(\psi_j)-\psi_j}_+}_{L^1} \le Cs.
    \end{align*}
    Since the density of $\psi_j$, $e^{\lambda_j\varphi_j}f_j$, has uniformly bounded $L^p$ norm by Lemma \ref{unifpotentialboundi},  Proposition \ref{GZProp} implies 
    \[
    W_s-\psi_j\le Cs^\gamma\le C\eta.
    \]

    Since $\psi_j$ is locally bounded for $x\in U$ and both 
    \[
    \Psi(\psi_j)\to \psi_j, \qquad -c_s\log\pare{\frac{t}{s}}\to \infty
    \]
    as $t\to 0$, we see that the infimum in the definition of $U_{c_s,s}$ is attained at some $t>0$. The lower bound in \eqref{curvest} implies that, pointwise for each $x\in X$, we have
    \[
    W_s(x)-\psi_j(x) \ge \pare{1-\eta}c_s\log\left(\frac{s}{t}\right)-Cs-\eta H(x).
    \]
    For this specific value $t$, we thus find constants $C_1,B_1\ge 0$ such that
    \[
     sC_1e^{-B_1H(x)}\leq t.
    \]
    Let $\xi_s=sC_1e^{-B_1H(x)}$. Using the monotonicity of the quantity $\Psi(\psi_j)\pare{x,w}+Kw^2$
    in $w$, together with \eqref{curvest} and \eqref{diNezzaTrapaniB}, we see
    \begin{align*}
    \Psi(\psi_j)\pare{x,\xi_s}-\psi_j(x) &\le \Psi(\psi_j)\pare{x,t}+Kt^2-\psi_j(x) \le U_{c_s,s}(x)+Cs-\psi_j(x).
    \end{align*}
    By the definition of $W_s$ given in \eqref{convexcombo}, we thus have 
    \begin{align*}
      \Psi(\psi_j)\pare{x,\xi_s}-\psi_j(x)  &\leq \frac{1}{1-\eta} \big(W_s(x)-\psi_j(x)+\eta(\psi_j(x)-\phi(x))\big)+Cs \\
        &\le \frac{\eta}{1-\eta}\pare{C+H(x)}+Cs \\
        &\le Cs^{\beta_0}\pare{1+H(x)}.
    \end{align*}
    Letting $r=\xi_s$, we then have
    \[
    \abs{\Psi(\psi_j)(x,w)-\psi_j(x)}\le Ce^{BH(x)}r^{\beta_0}.\qedhere
    \]
\end{proof}

\begin{remark}
    The preceding argument follows \cite{DemaillyDinewGuedjHiepKolodziejZeriahi}. As observed in \cite[Remark 3.2]{DiNezzaTrapani}, the Kiselman-Legendre transform used there, following \cite[Lemma 1.12]{BermanDemailly}, omits the linear correction term \(K(t-\delta)\). The corresponding estimates were corrected in \cite{KolodziejNguyen} by including this linear term, as we have done above. 
\end{remark}

\subsection{Preliminary Estimates along the Exhaustion}

\begin{cor}\label{uniformholderbound}
    There exist $0<\alpha<1$ and $B,C>0$ such that, for every sufficiently large exhaustion index $J$ 
    \[
    \sup_j\norm{\varphi_j}_{C^{0,\alpha}\pare{\overline{\Omega}_J}}\le Ce^{BJ}.
    \]
\end{cor}

\begin{proof}
    Let $J$ be chosen as in Lemma \ref{exhaustioncovering}. For $x\in \Omega_{J+1}$, the difference $H(x)=V_\theta(x)-\phi(x)$ satisfies the bound: 
    \[0\leq H(x)\leq C(J+1).\] 
    On the ball $B_{3R_J}\big(z_\ell^{(J)}\big)\Subset \Omega_{J+1}$, we may apply Lemma \ref{Holdercont} to get 
    \[
    |\Psi(\varphi_j)(x,\delta)-\varphi_j(x)|\leq Ce^{BJ}
    \delta^{\beta_0}
    \]
    for $0<\delta\leq R_J$.
    Now, by Lemma \ref{unifpotentialboundi} on $B_{3R_J}\big(z^{(J)}_\ell\big)$ we have 
    \[
    \|\varphi_j\|_{L^\infty(\Omega_{J+1})}\leq C(J+1).
    \]
    Hence, we see
    \[
    |D\Psi(\varphi_j)(z,\delta)|\leq \frac{Ce^{BJ}}{\delta}
    \]
    on $B_{2R_J}\big(z_\ell^{\pare{J}}\big)$.
    Let $x,y\in B_{2R_J}\big(z_\ell^{(J)}\big)$, and set \[\alpha=\frac{\beta_0}{1+\beta_0}.\] 
    We split into two cases:
    \begin{itemize}
        \item \textbf{Case 1:} If $|x-y|\le R_J^{1+\beta_0}$, let $\delta= |x-y|^{\frac{1}{1+{\beta_0}}}$. We can apply the local estimates above to get
    \begin{align*}
        \abs{\varphi_j(x)-\varphi_j\pare{y}} &
        \le 2Ce^{BJ}
        \delta^{\beta_0} + Ce^{BJ}
        \delta^{\beta_0} \leq Ce^{BJ}
        |x-y|^\alpha.
    \end{align*}
    \item \textbf{Case 2:} If $|x-y|>R_J^{1+{\beta_0}}$, we see 
    \begin{align*}
        \abs{\varphi_j(x)-\varphi_j\pare{y}} &\le 2\norm{\varphi_j}_{L^\infty(\Omega_{J+1})} \leq C(J+1)R^{-{\beta_0}}_J|x-y|^\alpha \leq Ce^{BJ}
        |x-y|^\alpha.
    \end{align*}
    \end{itemize}
    By the covering provided in Lemma \ref{exhaustioncovering} and the uniform estimates above, these local H\"older estimates extend to arbitrary points $x,y\in\overline{\Omega}_J$. To compare two nearby points, use the estimates on the doubled ball of a cover ball containing the first point; for points separated by at least a fixed multiple of the cover radius, use the supremum bound.
\end{proof}

\begin{remark}\label{limitworks}
    The upshot of these computations is that the family $\{\varphi_j\}_{j}$ is precompact in $L^1\pare{X}$ by Lemma \ref{unifpotentialboundi}. After passing to a subsequence, we have $\varphi_j\to\varphi_\infty \in\mathrm{PSH}(X,\theta)$ in $L^1\pare{X}$.

    Corollary \ref{uniformholderbound} gives equicontinuity and uniform boundedness on each fixed exhaustion level. By diagonal extraction, the subsequence converges locally uniformly on $U$, with the local limit equal to the restriction of the $L^1\pare{X}$ limit $\varphi_\infty$.
    
    Since $\lambda_j\searrow 0$ and the $\varphi_j$ are bounded in $L^1$, we have that ${\lambda_j}\varphi_j\to0$ in measure. Their uniform upper bound gives $e^{\lambda_j \varphi_j}$ a uniform upper bound as well. Thus, $e^{\lambda_j\varphi_j}\to 1$ in $L^{p^\prime}\pare{X,\omega^n}$, where $p^\prime=\frac{p}{p-1}$. H\"older's inequality therefore provides
    \[
    \norm{e^{\lambda_j\varphi_j}f_j-f}_{L^1\pare{X,\omega^n}}\le \norm{e^{\lambda_j\varphi_j}-1}_{L^{p^\prime}\pare{X,\omega^n}}\norm{f_j}_{L^{p}\pare{X,\omega^n}} + \norm{f_j-f}_{L^1\pare{X,\omega^n}}\to 0.
    \]
    On $U$, add a local smooth potential of $\theta$ and apply Bedford-Taylor continuity under local uniform convergence \cite{BT}. This gives $\inner{\theta^n_{\varphi_\infty}}=e^F\omega^n$ on $U$. Since neither measure puts mass on $Z$, this equality holds globally; see \cite{BEGZ} for further details. By uniqueness from \cite{BEGZ}, $\varphi$ and $\varphi_\infty$ differ by a constant. Throughout the remainder of the argument, we work with this subsequence.
\end{remark}

\section{Smoothing Strict Subsolutions}\label{section3}

The goal now is to obtain estimates for strict subsolutions that are compatible with the above exhaustion. The main difficulty is that, while one can obtain good estimates for a single ball via the work of \cite[Proposition 1.5]{EGZ2}, we need to patch the estimates carefully to ensure they are uniform over the entire cover. The Richberg technique, as employed in \cite{HLP}, provides a useful framework for moving these estimates from local to global calculations via regularization of the maximum function.

Let $K\Subset W\subset U$. By compactness, $K$ may be covered by $N$ holomorphic coordinate balls of radius $0<R\leq R_0$ such that
    \begin{equation}\label{specificcover}
        \overline{K}\subseteq
        \bigcup_{j=1}^{N}B_{R}(z_j),
        \qquad
        B_{3R}(z_{j})\Subset W.
    \end{equation}
For fixed $0<\alpha<1$, suppose $u\in C^{0,\alpha}(W)\cap \operatorname{PSH}(W,\theta)$ satisfies the Monge-Amp\`ere inequality
\begin{equation}\label{specificmaeqn}
\langle (\theta+\pot u)^n\rangle \geq f\omega^n.
\end{equation}
We previously considered $\psi\in C^\infty(W)$ satisfying $\theta_\psi\geq \delta \omega$.

\begin{lem}\label{C4subsolutions}
    On the cover given in \eqref{specificcover}, suppose that $u$ satisfies the inequality in \eqref{specificmaeqn} where $\inf_W f>0$ and $f^{1/n}\in C^{0,\alpha}(W)$. For $0<t,\eta\leq 1/4$, let 
    \[
    v=(1-t)u+t\psi.
    \]
    Then there are positive constants $c_*, c_1$ such that for every $0< b\leq \min\{1,c_*\eta t \delta R^2\}$ there exists a smooth function $p$ on a neighborhood of $K$ such that 
    \begin{enumerate}
        \item $v-2b\leq p\leq v$,
        \item $\theta+\pot p\geq c_1t\delta \omega$,
        \item $\det(\omega^{-1}(\theta+\pot p))^{1/n}\geq (1-\eta)(1-t)f^{1/n}$. 
    \end{enumerate}
    Here, $c_*,c_1$ depend only on $(X, \omega)$, $\theta$, and a chosen atlas of $X$.
\end{lem}

\begin{proof}

We begin with a local description on the $j$-th ball of the cover $B_R^{(j)}=B_R(z_j)$. We may represent the background metric $\omega$ via
\[\omega=\sqrt{-1}g_{k\bar{\ell}}\, dz_k\wedge d\bar{z}_\ell\]
and the closed form $\theta$ by a smooth potential $a_j$, i.e. on each ball $B_{3R}^{(j)}$ 
\[\theta=\sqrt{-1}\partial \overline \partial a_j.\] 
In particular, we see that $u+a_j$ is plurisubharmonic, and satisfies
\[
(\pot (u+a_j))^n\geq f \omega^n.
\]
Furthermore, the function $\psi+a_j$ is strictly plurisubharmonic. If we consider a non-negative hermitian matrix $Q$, a result of Gaveau \cite{Gaveau} gives
\[
\det(Q)^{1/n} = \frac1n\inf\{\operatorname{Tr}(QB): \det(B)=1, B\text{ is positive hermitian}\}.
\]
For every positive hermitian matrix $A$ with $\det(A)=1$, we then have
\begin{align*}
        \mathrm{Tr}(A(v+a_j)_{k\bar{\ell}})&= \mathrm{Tr}(A((1-t)(u+a_j)+t(\psi+a_j))_{k\bar{\ell}})\\
        &\geq (1-t)\mathrm{Tr}(A(u+a_j)_{k\bar{\ell}})\\
        &\geq n(1-t)f^{1/n}(\det(g_{k\bar{\ell}}))^{1/n}.
\end{align*}
Note that the above trace inequalities are taken in a distributional sense; see \cite[Proposition 1.5]{EGZ2} for further details. 

We now begin regularizing our subsolutions in $B_{2R}^{(j)}$. Let $\rho\in C_c^\infty(\mathbb{R})$ be a nonnegative, even mollifier satisfying $\int_{\mathbb{R}} \rho(x)\, dx=1$, and $\mathrm{Supp}(\rho)\subseteq [-1,1]$. We can then define for each $r>0$ 
\[\rho_r:\mathbb{R}^{2n}\to \mathbb{R}, \qquad \rho_r(x) =\frac{c_n}{r^{2n}}\rho\left(\frac{|x|}{r}\right)\]
where $c_n$ is a constant chosen to ensure $\int_{\mathbb{R}^{2n}}\rho_r(x)\, dx=1$, and we identify $\mathbb{C}^n\cong \mathbb{R}^{2n}$. For $0<r<R/10$, consider
\[
\hat v_j = ((v+a_j)* \rho_r)-a_j.
\]
We see
\[
\theta+\sqrt{-1}\partial \overline\partial  \hat{v}_j\geq c_0 t\delta \omega
\]
for some $c_0>0$. For small $r$, we then get the uniform bounds 
    \[
    \|\hat{v}_j-v\|_{L^\infty(B_{2R}^{(j)})}\leq b/10
    \]    
and the H\"older control
\[[f^{1/n}(\det(g_{k\bar{\ell}}))^{1/n}]_{C^{0,\alpha}(B_{3R}^{(j)})}r^\alpha\leq \frac{\eta}{4}\inf_{B_{3R}^{(j)}}\left\{
f^{1/n}(\det(g_{k\bar{\ell}}))^{1/n}\right\}. \]
Hence,
\begin{equation}\label{vhatmetricbound}
\det\left(
        \omega^{-1}\left(
            \theta+i\partial\bar\partial\hat{v}_j
        \right)
    \right)^{1/n}
    \geq (1-\eta/4)(1-t)f^{1/n}.
\end{equation}

We now want to make these estimates amenable to the given cover of $K$. To this end, consider the smooth cutoff functions $\chi_j\in C^\infty(B_{2R}^{(j)})$ satisfying $|\chi_j|\leq 3$ and
    \[
    \chi_j(z) =\begin{cases}
        3, &z\in B_R^{(j)},\\
        -3,& z\in B_{2R}^{(j)}\setminus B_{3R/2}^{(j)}\\
    \end{cases}
    \]
    along with the derivative bounds
    \[|D^k\chi_j|\leq C_kR^{-k}\] 
    for $1\leq k\leq 4$. In particular, we have
    \[
    |\pot \chi_j |\leq \frac{C}{R^2}.
    \]
    If we  define
    \[
    \kappa_j = \hat{v}_j+b\chi_j
    \]
    then we see 
    \begin{equation}\label{kappajmetricbound}
    \theta+\sqrt{-1}\partial\overline\partial \kappa_j \geq (1-\eta/4)(\theta+\sqrt{-1}\partial\overline\partial \hat{v}_j).
    \end{equation}
    Since $(1-\eta/4)^2\geq 1-\eta$, we can combine \eqref{vhatmetricbound} and \eqref{kappajmetricbound} to obtain
    \[
    \det(\omega^{-1} (\theta+i\partial \overline\partial \kappa_j))^{1/n}\geq (1-\eta)(1-t)f^{1/n}.
    \]

    To finish, we  wish to obtain uniform control over our cover by considering how to choose the maximum contribution across all balls $B_{2R}^{(j)}$. Let $\tau_j\in C^\infty(B_{2R}^{(j)})$ be another smooth cutoff function supported in $B_{15R/8}^{(j)}$, satisfying $0\leq \tau_j\leq 1$ and $\tau_j(z)=1$ for $z\in B_{7R/4}^{(j)}$. If $D=\bigcup_j \overline{B}_{R}^{(j)}$, then 
    \[
    L = \inf_D v-4b.
    \]
    Define now
    \[
    w_j(z) = \begin{cases}
        L+\tau_j(\kappa_j(z)-L), &z\in B_{2R}^{(j)}\\
       L, &z\in D\setminus B_{2R}^{(j)}
    \end{cases}
    \]
    This choice of $L$ ensures that $L$ remains a uniform lower bound for $w_j$. 
    Also, on $B_R^{(j)}$
    \begin{equation}\label{lowerboundwj}
        v+2.9b\leq w_j
    \end{equation}
    while on $B_{2R}^{(j)}\setminus B_{7R/4}^{(j)}$ we have 
    \begin{equation}\label{upperboundwj}
    w_j\leq v-2.9b. 
    \end{equation}
    Using the aforementioned mollifier $\rho$, define the regularized maximum $M:\mathbb{R}^N\to \mathbb{R}$ as follows:
    \[
    M(y_1, \dots, y_N) = \int_{\mathbb{R}^N} \max_{1\leq j\leq N}\left(y_j+\frac{b}{10}\zeta_j\right)\prod_{j=1}^N\rho(\zeta_j)\, d\zeta.
    \]
    As in \cite{HLP}, $M$ is smooth and convex. By construction, we have
    \[
    |M(y)-\max_{1\leq j\leq N} y_j|\leq \frac{b}{10}. 
    \]
    Furthermore, since $M$ is monotone increasing in each $y_j$, we see 
    \[
    \frac{\partial M}{\partial y_j}\geq 0, \qquad 
    \sum_{j=1}^N \frac{\partial M}{\partial y_j}=1.
    \]
    We may then define on the cover $D$
    \[p=M(w_1, \dots, w_N)-4b.\]
    Note that 
    \[
    v+2.8b\leq M(w_1, \dots, w_N) \leq v+3.2b.
    \]
    Hence, $p$ satisfies condition (1). 
    
    For each $x\in D$, if $x\notin B_{7R/4}^{(j)}$ for some $j$, then by \eqref{lowerboundwj} and \eqref{upperboundwj}, there is an index $\ell$ such that 
    \[
    w_j(x)+\frac{b}{10}\zeta_j\leq v(x)-2.8b<v(x)+2.8b\leq w_\ell(x)+\frac{b}{10}\zeta_\ell
    \]
    for every $|\zeta_j|,|\zeta_\ell|\leq 1$, and we consider constant extensions of $w_j$ and $w_\ell$ as necessary. 
    Thus, the $j$-th argument of $M(w_1, \dots, w_N)$ does not contribute to $M$. Hence, the only way $w_j(x)$ contributes is if $x\in B_{7R/4}^{(j)}$, in which case it agrees locally with $\kappa_j$.
    By the chain rule,
    \begin{align*}
    \theta_p &= \sum_{j=1}^N \frac{\partial M}{\partial y_j}(w_1, \dots, w_N)\theta_{w_j} +\sqrt{-1}\sum_{j,k=1}^{N} \left(\frac{\partial^2 M}{\partial y_j\partial y_k}(w_1, \dots, w_N)\right)\partial w_j\wedge \overline \partial w_k.
    \end{align*}
    By the concavity and monotonicity of the operator $\Gamma(A) = \det(A)^{1/n}$ for positive hermitian matrices $A$, we deduce (2) and (3). 
    \end{proof}

\begin{cor}\label{C4bound}
    For each fixed $\alpha$, the smooth function $p$ given by Lemma \ref{C4subsolutions} obeys
    \[
    \|p\|_{C^4(K)}\leq Q
    \]
    where $Q$ is a polynomial depending on 
    \[\left\{\alpha,\quad N,\quad\|v\|_{C^{0, \alpha}(W)},\quad \|f^{1/n}\|_{C^{0,\alpha}(W)},\quad \frac1R,\quad \frac1b,\quad \frac1t,\quad \frac1\eta,\quad \frac1\delta,\quad \frac{1}{\inf_W f}\right\}. 
    \] 
\end{cor}

\begin{proof}
Firstly, note that $r$ can be chosen to be a sufficiently small fixed multiple of 
\[
\min\left\{1,R,\left(\frac{b}{1+\|v\|_{C^{0,\alpha}(W)}}\right)^{1/\alpha}, \left(\frac{\eta\min\{1, (\inf_W f)^{1/n}\}}{1+\|f^{1/n}\|_{C^{0,\alpha}(W)}}\right)^{1/\alpha}\right\}.
\]
The convolution $\hat{v}_j$ obeys
\[
\|\hat{v}_j\|_{C^4(B_{2R}^{(j)})}\leq A\left(1+\frac1{r^4}\right)
\]
where $A=A(1+\|v\|_{L^\infty (W)})$, while the cutoff $\kappa_j$ obeys 
\begin{equation}\label{mess}
\|\kappa_j\|_{C^4(B^{(j)}_{2R})}\leq \left(1+\frac{1}{R^4}\right)\left[A\left(1+\frac{1}{r^4}\right)+b\left(1+\frac{1}{R^4}\right)\right].
\end{equation}
By the construction of $M$, we see that
\[
0\leq\frac{\partial M}{\partial y_j}\leq 1.
\]
Let $I$ be a multiindex with $1\leq |I|\leq4$. We can show that 
\begin{equation}\label{maxestimate}
\left|\frac{\partial^{|I|} M}{\partial y_I}\right|\leq C_I \left(\frac{b}{10}\right)^{1-|I|}
\end{equation}
by putting one weak derivative on the maximum inside $M$, and the other derivatives on the mollifiers $\rho$. Only the arguments $w_j$ of $M$ which agree locally with $\kappa_j$ contribute to $M$, and are hence controlled by \eqref{mess}. By repeatedly applying the chain rule, if we denote the upper bound in \eqref{mess} by $D=D(1/R, 1/r, A, b)$, we see by \eqref{maxestimate}
\[
\|p\|_{C^4(K)}\leq A+C(1+N)^4\left(1+\frac{1}{b^3}\right)(1+D)^4
\]
where $C$ is independent of $N$. Substituting the displayed choice of $r$ into these estimates and rounding the exponent $1/\alpha$ to the integer $\lceil1/\alpha\rceil$ gives the claimed polynomial bound. 
\end{proof}

On the exhaustion, set $K=\overline{\Omega}_{j}$, $W=\Omega_{j+1}$ and $t_j=e^{-3j}=\eta_j$. Define the constants
\[
a_j=-\frac{2n}{\lambda_j}\log\pare{1-t_j},\qquad b_j=ct^2_j\delta R_j^2
\]
for some constant $0<c<c_*$ from Lemma \ref{C4subsolutions}.
In particular, we have $a_j\to 0$ as $j\to \infty$.

If in the setting of Lemma \ref{C4subsolutions} we take $u=\varphi_j$, $v_j=\pare{1-t_j}\varphi_j+t_j\psi$, and $f=e^{\lambda_j\varphi_j+F_j}$, we have on $W$,
\[
0\le \varphi_j-v_j\le C\pare{1+j}t_j,\quad \theta_{v_j}\ge t_j\delta \omega,\quad \theta_{v_j}^n\ge \pare{1-t_j}^ne^{\lambda_j\varphi_j+F_j}\omega^n.
\]
The determinant inequality follows by retaining the $\theta_{\varphi_j}^n$ term in the mixed Monge–Ampère expansion. Lemma \ref{Fapprox} and the bound $|\varphi_j|\le C(1+j)$ on $W$ give $e^{\lambda_j\varphi_j+F_j}\ge ce^{-j}$. Corollary \ref{uniformholderbound}, the bounds for $F_j$ and $\nabla F_j$, and the mean value theorem for the exponential give $\|\pare{e^{\lambda_j\varphi_j+F_j}}^{1/n}\|_{C^{0,\alpha}(W)}\le Ce^{Bj}$. The same exponential bound holds for the \(C^{0,\alpha}\) norm of \(v_j\), using the barrier derivative bounds on the next exhaustion level. Thus Lemma \ref{C4subsolutions} and Corollary \ref{C4bound} apply with all the data and required reciprocals bounded exponentially in $j$. They give a smooth $p_j$ near $K$ satisfying
\[
v_j-2b_j\le p_j\le v_j,\quad \theta_{p_j}\ge c_1t_j\delta\omega, \quad \theta_{p_j}^n\ge \pare{1-t_j}^{2n}e^{\lambda_j\varphi_j+F_j}\omega^n,\quad \norm{p_j}_{C^4\pare{K}}\le Ce^{Bj}.
\]
Consider the shifted regularization
\begin{equation}\label{wjdefn}
    w_j = p_j-a_j.
\end{equation}
Since $p_j\le \varphi_j$ and $e^{-\lambda_j a_j}=\pare{1-t_j}^{2n}$, we see that $\theta_{w_j}^n\ge e^{\lambda_j w_j+F_j}\omega^n$.

We can rephrase the results of this section more conveniently as follows:

\begin{cor}\label{wpolynomialbound}
    Each $w_j$ as defined in \eqref{wjdefn} is smooth on a neighborhood of $\overline{\Omega}_j$, and there are fixed constants $c,C>0$ such that on $\overline{\Omega}_j$
    \begin{enumerate}
        \item $\varphi_j-\epsilon_j\leq w_j\leq \varphi_j$, where $\epsilon_j=Ce^{-3j}\pare{j+1}+2c t_j^2 \delta R_j^2+a_j$
        \item $\theta_{w_j}\geq ce^{-Cj}\omega$
        \item $\theta_{w_j}^{n}\geq e^{\lambda_jw_j+F_j}\omega^n$
        \item $\|w_j\|_{C^4(\overline{\Omega}_j)}\leq Ce^{Bj}
        $.
    \end{enumerate} 
\end{cor}

\begin{remark}
    The quantity $\epsilon_j$ records the errors at each stage of the regularization above. Importantly, each error decays exponentially as $j\to \infty$.
\end{remark}

\section{The Dirichlet Problem and Gradient Estimates}\label{section4}

The goal of this section is to leverage the analytic estimates of the $w_j$ on each exhaustion domain $\Omega_j$ to establish first-order estimates for the solutions to the Dirichlet problems of the form
\begin{equation}\label{DirichletProblem}
\begin{cases}
    \langle (\theta+i\partial\overline{\partial} u_j)^n\rangle =e^{\lambda_ju_j+F_j}\omega^n,&\text{on }\Omega_j\\
    u_j=w_j,&\text{on }\partial \Omega_j
\end{cases}
\end{equation}

\subsection{Boundary Geometry of the Exhaustion Domains}

\begin{lem}\label{Smoothexhaustion}
For sufficiently large integers $j$, each $\del\Omega_j$ is smooth.
\end{lem}

\begin{proof}
    Fix a log-resolution $\pi\colon \hat{X}\to X$ on $(\hat{X},\hat{\omega})$ as in the proof of Lemma \ref{psibounded}. On the band $B_j=\{x\in X:\abs{\psi+j}\le 2\}$, we have
    \begin{equation*}
        \sum_{\ell=1}^r (1+\sigma)a_\ell\log\abs{z_\ell}^2 = -j+O\pare{1}.
    \end{equation*}
    Because there are finitely many coordinates and charts,  for $j$ sufficiently large, at every $x\in B_j$, at least one coordinate of the resolution satisfies
    \begin{equation*}
        \abs{z_\ell}\le Ce^{-Bj}.
    \end{equation*}
    In that coordinate,
    \begin{equation*}
        \frac{\del \pi^*\psi}{\del z_\ell} = \frac{(1+\sigma)a_\ell}{z_\ell}+\frac{\del\tilde{h}}{\del z_\ell},
    \end{equation*}
    where $\tilde h=(1+\sigma) h$.
    The first term grows exponentially while the second term is bounded on each chart. Therefore, in sufficiently deep bands
    \begin{equation*}
        \abs{\frac{\del \pi^*\psi}{\del z_\ell}}_{\hat{\omega}} \ge 1.
    \end{equation*}
    Since $d\pi$ is bounded on the compact manifold $\hat{X}$, the chain rule yields 
    \begin{equation*}
        \abs{d\pi^*\psi}_{\hat{\omega}} \le C\abs{d\psi}_\omega.
    \end{equation*}
    Combining these results, on $B_j$ we see that there is a fixed $c_0>0$ such that 
    \begin{equation}
        \abs{\nabla\psi}_\omega \ge c_0
    \end{equation}
    for all sufficiently large $j$. Hence, $-j$ is a regular value of $\psi$ and 
    \[\del\Omega_j=\brac{x\in X:\psi(x)=-j}\] 
    is smooth as claimed.
\end{proof}

\begin{lem}[]\label{backgroundgeo}
For sufficiently large integers $j$ and for each positive integer $k$, the boundary charts of $\Omega_j$ and their signed-distance functions in suitable tubular neighborhoods have fixed-order bounds at most $C_ke^{B_kj}$. The reciprocals of the chart radii and tubular widths satisfy bounds of the same form.
\end{lem}

\begin{proof} We split the proof into three parts.

    \textit{Boundary chart:} Use a fixed finite collection of real coordinate charts for the background metric with uniform bounds for the metric and its derivatives.
    Set $f_j=\psi+j$ and choose $1\le M_j\le N_j
    $ bounding the first two derivatives of $f_j$ in the level band $\brac{\abs{f_j}\le 2}$. Coordinate versions of $\abs{\nabla\psi}_\omega$ have a fixed lower bound $c_0$, after perhaps changing fixed constants.
    
    Fix $p\in \del\Omega_{j}$, and translate and rotate coordinates so that $p=0$ and that the last coordinate direction agrees with the gradient at $p$. In these coordinates, at $\pare{x^\prime,x_{2n}}:=\pare{x_1,\ldots,x_{2n-1},x_{2n}}$
    \begin{equation*}
        f_j\pare{0}=0, \qquad D_{x^\prime}f_j\pare{0}=0, \qquad \del_{x_{2n}}f_j\pare{0} \ge c_0.
    \end{equation*}
    Choose $r_j=cM_j^{-1}$ with $c>0$ small enough that the cylinder
    \begin{equation*}
        \abs{x^\prime}<2r_j,\qquad \abs{x_{2n}}<2r_j
    \end{equation*}
    lies in the band. Such a $c$ exists because if the cylinder were to leave the band, we would have both $\abs{f_j}=2$ and by gradient bounds, $\abs{f_j}\le CM_jr_j$. Taking $c$ sufficiently small forces the second quantity to be less than $2$. 
    
    The Hessian bound gives 
    \begin{equation*}
        \del_{x_{2n}}f_j \ge \frac{c_0}{2},
    \end{equation*}
    throughout the cylinder, after perhaps shrinking $c$. Hence, every vertical line in the cylinder has at most one zero. Taylor's formula along the tangential plane gives
    \begin{equation*}
        \abs{f_j\pare{x^\prime,0}} \le CM_j\abs{x^\prime}^2.
    \end{equation*}
    For $\abs{x^\prime}<r_j$, if we shrink $c$ further, we can achieve
    \begin{equation*}
        \abs{f_j\pare{x^\prime,0}} \le \frac{c_0r_j}{4}.
    \end{equation*}
    The vertical derivative bound then implies
    \begin{equation*}
        f_j\pare{x^\prime,-r_j} < 0 < f_j\pare{x^\prime,r_j}.
    \end{equation*}
    Thus, $f_j\pare{x^\prime,\cdot}$ has exactly one zero, $x_{2n}=h_p\pare{x^\prime}$, and the implicit function theorem tells us that $h_p$ is smooth. Because this is the only zero, the level set in this inner cylinder is the boundary graph. Differentiating $f_j\pare{x^\prime,h_p\pare{x^\prime}}=0$ gives
    \begin{equation*}
        \del_i h_p = -\frac{\del_i f_j}{\del_{x_{2n}}f_j}.
    \end{equation*}
    Higher derivatives of $h_p$ can be given in terms of the derivatives of $\psi$ and powers of $\pare{\del_{x_{2n}}\psi}^{-1}$, which all have exponential bounds in $j$. So, for each fixed $k\ge 1$,
    \begin{equation*}
        \norm{h_p}_{C^k} \le C_ke^{B_kj}.
    \end{equation*}

    \textit{Distance function:} Consider the normal exponential map 
    \begin{equation*}
        \mathcal{N}_p\pare{x^\prime,t} = \exp_{\pare{x^\prime,h_p\pare{x^\prime}}}\pare{t\nu_j\pare{x^\prime,h_p\pare{x^\prime}}},
    \end{equation*}
    where $\nu_j$ is the inward unit normal. At $\pare{0,0}$, its differential maps the $x^\prime$ directions to tangent directions and the $t$ direction to the unit normal. Fixed metric equivalence gives a uniform positive lower bound $a$ for its smallest singular value.

    The graph estimates give $\abs{D\nu_j}\le CM_j$ and the exponential map has fixed derivative bounds near its zero section. Thus, for $\abs{x^\prime}<cr_j$ and $\abs{t}<cr_j$, we have
    \begin{equation*}
        \norm{D\mathcal{N}_p\pare{x^\prime,t}-D\mathcal{N}_p\pare{0,0}} \le \frac{a}{2}.
    \end{equation*}
    In ambient coordinates the parameter cylinder is convex. For two parameter points $z,z^\prime$ in the cylinder, integrating the derivative along their joining segment gives
    \begin{equation*}
        \abs{\mathcal{N}_p(z)-\mathcal{N}_p\pare{z^\prime}} \ge \frac{a}{2}\abs{z-z^\prime}.
    \end{equation*}
    Thus, $\mathcal{N}_p$ is injective on an entire uniform parameter cylinder and its differential has a uniformly bounded inverse.

    To get global injectivity, choose $\eta_j=cr_j$. Suppose there are two points $\pare{p,t}$ and $\pare{q,s}$ such that
    \begin{equation*}
        \exp_p\pare{t\nu_j\pare{p}}=\exp_q\pare{s\nu_j\pare{q}}, \qquad \abs{t},\abs{s}<\eta_j.
    \end{equation*}
    The lengths of the normal segments give $d_\omega\pare{p,q}\le \abs{t}+\abs{s}<2\eta_j$. So, $q$ lies in the single graph chart centered at $p$. Injectivity in the cylinder then gives $p=q$ and $t=s$. Hence
    \begin{equation*}
        \mathcal{N}_j\colon \del\Omega_j\times \pare{-\eta_j,\eta_j}\to X
    \end{equation*}
    is globally injective. 

    Now take $\eta_j$ below the fixed ambient injectivity radius. If $x=\mathcal{N}_j\pare{p,t}$, then 
    \begin{equation*}
        \operatorname{dist}_\omega\pare{x,\del\Omega_j}\le \abs{t},
    \end{equation*}
    because the normal segment is a path from $p$ to $x$. By the compactness of $\del\Omega_j$, there is a closest point $q\in \del \Omega_j$. A minimizing geodesic from $q$ to $x$ is perpendicular to $\del\Omega_j$ and thus,
    \begin{equation*}
        x = N_j\pare{q,s},\qquad \abs{s}=\operatorname{dist}_\omega\pare{x,\del\Omega_j}\le \abs{t}.
    \end{equation*}
    By global injectivity, $p=q$, $t=s$, and $\operatorname{dist}_\omega\pare{N_j\pare{p,t},\del\Omega_j}=\abs{t}$. Conversely, because the unit normal points inward, the same nearest-point argument represents every point of distance less than $\eta_j$ in this way.
    Thus the signed distance, positive inside $\Omega_j$, satisfies
    \begin{equation*}
        d_j\pare{N_j\pare{p,t}} = t.
    \end{equation*}
    Consequently, the signed-distance function satisfies
    \begin{equation*}
        \norm{d_j}_{C^k} \le C_ke^{B_kj}.
    \end{equation*}
    
    \textit{Chart radii:} We have that $\eta_j^{-1}=Cr_j^{-1}\le N_j$.
\end{proof}

\subsection{Existence Results for the Dirichlet Problem}

\begin{prop}\label{ujexist}
    For sufficiently large $j$, \eqref{DirichletProblem} has an admissible smooth solution $u_j\in C^\infty(\overline{\Omega}_j)$.
\end{prop}

\begin{proof}
    By Lemma \ref{Smoothexhaustion}, for sufficiently large $j$, the boundary $\del\Omega_j$ is smooth. By Corollary \ref{wpolynomialbound}, we see that $w_j$ is a strictly admissible subsolution to the equation 
    \[
    \pare{\theta+\pot u}^n = e^{\lambda_ju+F_j}\omega^n
    \]
    with $u(z) = w_j(z)$ for $z\in \partial \Omega_j$. Additionally, $w_j$ is smooth on $\overline{\Omega}_j$. The right-hand side is smooth, positive, and non-decreasing in $t$. The smooth admissible subsolution $w_j$ has the prescribed boundary trace and satisfies the determinant inequality on \(\overline{\Omega_j}\). Thus \cite[Theorem 1.1]{GuanLi} applies, on each component if necessary, and supplies a solution $u_j$, smooth up to the boundary.
\end{proof}

The Dirichlet solutions $u_j$ are also close to the global potentials $\varphi_j$.

\begin{lem}\label{ujclose}
    For $j$ sufficiently large,
    \begin{equation}\label{inequalities}
    \varphi_j-\epsilon_j\le w_j\le u_j\le \varphi_j,\quad\text{on }\overline{\Omega}_j.
    \end{equation}
\end{lem}
\begin{proof}
    
    The first inequality follows from Corollary \ref{wpolynomialbound}. For the second inequality, suppose that $w_j-u_j$ attains a positive maximum at an interior point $x_0\in \Omega_j$. Then $\theta_{w_j}(x_0)\le \theta_{u_j}(x_0)$ while the subsolution inequality gives
    \begin{align*}
        \theta_{w_j}^n \ge e^{\lambda w_j(x)+F_j(x)}\omega^n > e^{\lambda u_j(x)+F_j(x)}\omega^n = \theta_{u_j}^n
    \end{align*}
    contradicting the above inequality. Hence $w_j\le u_j$.

    For the last inequality, note that $u_j\le \varphi_j$ on the boundary. Both $u_j$ and $\varphi_j$ are bounded and continuous on $\overline{\Omega}_j$, so Bedford-Taylor comparison \cite{BT} gives
    \[
    \int_{\brac{\varphi_j<u_j}}\theta_{u_j}^n\le \int_{\brac{\varphi_j<u_j}}\theta_{\varphi_j}^n.
    \]
    Assuming the set of integration is nonempty, it would be open, and the equations with $\lambda_j>0$ would imply that the left integral is larger than the right integral. This is impossible, so $u_j\le \varphi_j$.
\end{proof}

Lemma \ref{unifpotentialboundi} and \eqref{inequalities} combine to give the following corollary:

\begin{cor}\label{detuniformity}
    Define $\kappa=\frac{\sigma}{1+\sigma}$. Then there is a uniform constant $C>0$ such that for each $j$,
    \begin{equation}
        \phi-C\le u_j\le C,\qquad v_j:=u_j-\psi\ge -\kappa \psi-C.
    \end{equation}
\end{cor}

\subsection{Gradient Control of the Dirichlet Solutions}

We now use the regularity of $w_j$ to establish gradient estimates for $u_j$.

\begin{lem}\label{BoundaryC1Estimates}
    There exist uniform constants $C,L_1>0$ such that for sufficiently large $j$, \[\sup_{\partial\Omega_j}\abs{\nabla u_j}^2_\omega\le Ce^{L_1j}.\]
\end{lem}
\begin{proof}
    Set $z_j=u_j-w_j$, so that $z_j(p)=0$ for all $p\in \partial \Omega_j$. Since $u_j$ is admissible,
    \begin{align*}
        \Delta_\omega z_j &= \operatorname{Tr}_\omega\theta_{u_j} - \operatorname{Tr}_\omega\theta_{w_j} \ge  - \operatorname{Tr}_\omega\theta_{w_j} \ge -Ce^{Bj}=-N_j.
    \end{align*}
    In particular, the final inequality follows from the $C^4$-bound from Corollary \ref{wpolynomialbound}. For notational convenience,
    let $d_j$ denote the inward distance to $\partial\Omega_j$ on a tubular neighborhood. After perhaps increasing $C$, choose $N_j$ sufficiently large so that
    \begin{equation}
       0\leq d_j\leq \frac1{N_j}, \qquad  0\le z_j\le N_j,
    \end{equation}
    and the distance obeys
    \begin{equation}\label{distancebounds}
        \abs{\Delta_\omega d_j}\le N_j, \qquad \abs{\nabla d_j}^2_\omega\ge c_0>0
    \end{equation} 
    by Lemma \ref{backgroundgeo} for some constant $c_0>0$. For constants $c,C_1>0$, we  construct a barrier $b_j$ by defining the following:
    \begin{align*}
        \eta_j = \frac{c}{N_j^2},\qquad A_j=\frac{4N_j}{\eta_j},\qquad B_j= C_1(N_j+A_jN_j), \qquad b_j(d) = A_jd-B_jd^2.
    \end{align*}

    Choose $C_1$ so that $2c_0C_1\ge 1$, and choose $c$ sufficiently small that $\eta_j$ lies in the distance collar and the inequalities
    \[c\le \frac{1}{2C_1}\qquad B_j\eta_j\le \frac{A_j}{2}\] 
    hold. Note that by our control on the Laplacian and distance collar these choices are independent of $j$.
    
    On the inner collar $0\le d_j\le \eta_j$, we see that 
    \begin{equation}\label{Aest}
        0\le A_j-2B_jd\le A_j.
    \end{equation}
    Equations \eqref{distancebounds} and \eqref{Aest} combine to bound the following:
    \begin{align*}
        \Delta_\omega b_j &= \pare{A_j-2B_jd_j}\Delta_\omega d_j -2B_j\abs{\nabla d_j}^2_\omega \le A_jN_j-2c_0B_j \le -N_j.
    \end{align*}
    On the actual boundary $d_j=0$, we see that $b_j$ and $z_j$ both vanish. On the inner collar boundary $d_j=\eta_j$,
    \begin{align*}
        b_j\pare{\eta_j}=\eta_j\pare{A_j-B_j\eta_j}\ge \frac{A_j\eta_j}{2}=2N_j\ge z_j(\eta_j).
    \end{align*}
    Because $\Delta_\omega\pare{z_j-b_j}\ge 0$, the maximum principle implies that $0\le z_j\le b_j$ on the collar. Therefore the inward normal derivative of $z$ at the boundary lies in $\sqbr{0,A_j}$, while its tangential derivatives vanish. Thus,
    \begin{align*}
        \sup_{\partial \Omega_j}\abs{\nabla u_j}_\omega \le \norm{w_j}_{C^1}+A_j\le CN_j^3\leq C'e^{B' j}.
    \end{align*}
    Squaring this relation and adjusting the exponents gives the uniform constant $L_1$ in the statement.
\end{proof}

\begin{lem}\label{determinantbounded}
   There exist constants $c,C$ independent of $j$ such that 
   \[
    ce^{-Cj}\le \det\omega^{-1}\theta_{u_j}\le Ce^{Cj} \qquad \text{on }\overline{\Omega}_j.
    \]
\end{lem}
\begin{proof}
    We see that there is a uniform $C>0$ such that
    \[
    -\frac{j}{1+\sigma}-C\le u_j\le C
    \]
    on $\overline{\Omega}_j$. Because $\lambda_j=e^{-j}$, there is some constant $C_0$ uniform in $j$ such that $\abs{\lambda_j u_j}\le C_0$. The claim then follows from Lemma \ref{Fapprox}. 
\end{proof}

The comparison with $w_j$ and $\varphi_j$ gives uniform bounds for $u_j$ relative to the barrier and controls the boundary gradient. We now combine these facts with the exponential bound for $\nabla F_j$ to control the gradient throughout $\Omega_j$.

\begin{prop}\label{GradientEstimateuj}
    For sufficiently large $j$, there exist  uniform constants $C,C^\prime>0$ such that $\abs{\nabla u_j}_\omega^2\le Ce^{C^\prime j}$ on each $\Omega_j$.
\end{prop}
\begin{proof}
    We follow a similar technique to \cite{PS2}. To simplify the notation, we  suppress the subscripts $j$. Let $m$ be a uniform lower bound for $v$ as defined in Corollary \ref{detuniformity}. If $S=\abs{\nabla u}^2_\omega$, define 
    \begin{align*}
        t = v-m+1\ge 1, \qquad h(v)=-Bv+t^{-1}, \qquad p=-h^\prime=B+t^{-2}, \qquad Q=e^{h(v)}S.
    \end{align*}
    It suffices to show that $Q\le Ce^{bj}$. Since $u$ has a uniform upper bound and $\psi\ge -j$ on $\overline{\Omega}_j$, we have $e^{-h}\le Ce^{Bj}$, and after adjusting the exponent we have the desired bound on $S$.

    Consider the point $x_0\in \Omega$ where $Q$ attains a maximum. If $x_0\in \partial \Omega$, by Lemma \ref{BoundaryC1Estimates}, we have  
    \[
    Q\le C_Be^{\pare{L_1-B\kappa}j}.
    \]
    Taking $B\ge \frac{L_1}{\kappa}$ gives the desired upper bound for $Q$. If $x_0$ lies in the interior and $S<1$, because $h$ has a uniform upper bound, $Q$ is also sufficiently bounded. We therefore assume that $x_0$ lies in the interior and $S\ge 1$ at $x_0$.
    
    Choose normal holomorphic coordinates for $\omega$ that diagonalize $\theta_u$ at $x_0$. 
    We can write the metric locally by
    \begin{equation}\label{metricdef}
    \tilde{g}_{i\bar{k}} =\theta_{i\bar{k}}+u_{i\bar{k}}.
    \end{equation}
    If $a_i=\tilde{g}^{i\bar{i}}$,  rewrite the Laplacian and trace in the diagonalization as 
    \begin{align*}
        L=\sum_i a_i\del_i\del_{\bar{i}}, \qquad T=\sum_i a_i.
    \end{align*}
    We also define a bookkeeping term
    \begin{align*}
        K=\sum_{i,k}a_i\pare{\abs{u_{ki}}^2+\abs{u_{k\bar{i}}}^2}.
    \end{align*}
    At the interior point $x_0$, we have $LQ\le 0$. We now expand the expression $LQ$. 
    First examining $LS$, we see that
    \begin{align*}
        LS &= K + 2\operatorname{Re}\left(\sum_{i,k}a_iu_ku_{i\bar{i}\bar{k}}\right)+\sum_{i,k,\ell}a_i\pare{\del_i\del_{\bar{i}}g^{k\bar{\ell}}}u_ku_{\bar{\ell}}.
    \end{align*}
    Because the background metric is fixed and smooth, 
    \[
    \abs{\sum_{i,k,\ell}a_i\pare{\del_i\del_{\bar{i}}g^{k\bar{\ell}}}u_ku_{\bar{\ell}}}\le C\pare{\sum_ia_i}\pare{\sum_k\abs{u_k}^2}=CST.
    \]
    We would like to show a similar bound for the second term of $LS$. From $\det \tilde{g} = e^{\lambda u + F}\det g$, differentiating in the $\bar{k}$ direction yields
    \begin{align*}
        \sum_{i}a_iu_{i\bar{i}\bar{k}} &= \lambda u_{\bar{k}} + F_{\bar{k}} + \del_{\bar{k}}\log \det g - \sum_i a_i\theta_{i\bar{i}\bar{k}}.
    \end{align*}
    Thus, 
    \begin{align*}
        LS &\ge K + 2\operatorname{Re} \left(\sum_{k}u_k\pare{\lambda u_{\bar{k}} + F_{\bar{k}} + \del_{\bar{k}}\log \det g}\right) - 2\operatorname{Re}\left(\sum_{i,k} a_iu_{{k}}\theta_{i\bar{i}\bar{k}}\right) - CST.
    \end{align*}
    The first sum is bounded from below by $-C\pare{1+S}-Ce^{Cj}S$, and by Cauchy-Schwarz, the second sum is bounded below by $-C\pare{1+S}T$. Note that in the first sum, we are using the fact that $\norm{\nabla F}_{L^\infty\pare{\Omega,\omega^n}}\le Ce^{Cj}$ and that $\lambda S\ge 0$. 
    Because $S\ge 1$, we can say $LS\ge K-CST-Ce^{Cj}S$. Similarly,
    \begin{align*}
        Lv = n-\operatorname{Tr}_{\theta_u} \theta_\psi \le n-\delta T,\qquad e^{-h}Le^{h}\ge \pare{p^2+h^{\prime\prime}}\abs{\nabla v}^2_{\theta_u}-np+\delta pT.
    \end{align*}
    Combining everything, we see that
    \begin{equation}\label{estimate4.1}
        0\ge LS+S\pare{\pare{p^2+h^{\prime\prime}}\abs{\nabla v}^2_{\theta_u}-np+\delta pT} -2p\operatorname{Re}\left(\sum_i a_iv_iS_{\bar{i}}\right).
    \end{equation}
    At the center of the normal coordinates we see that 
    \begin{align*}
        S_{\bar{i}} &= \sum_{k}u_{k\bar{i}}u_{\bar{k}}+u_k\overline{u_{ki}}.
    \end{align*}
    Cauchy-Schwarz and Young's inequality bound the second part of the cross term below by 
    \begin{align*}
        -\sum_{i,k}a_i\abs{u_{ki}}^2-p^2S\abs{\nabla v}^2_{\theta_u},
    \end{align*}
    Locally, by \eqref{metricdef}, exploiting diagonality we get
    \begin{align*}
        \operatorname{Re}\left(\sum_{i,k} a_i v_iu_{k\bar{i}}u_{\bar{k}}\right) &=
        \operatorname{Re}\left(\sum_{i,k} a_i v_i(\tilde{g}_{k\bar{i}}-\theta_{k\bar{i}})u_{\bar{k}}\right)\\
        &=\operatorname{Re}\left(\sum_i v_i u_{\bar{i}}\right)-\operatorname{Re}\left(\sum_{i,k}a_i v_i\theta_{k\bar{i}}u_{\bar{k}}\right)\\
        &=(I)+(II).
    \end{align*}
   By definition of $v$, we have 
    \begin{equation}\label{quant1}
    (I)\leq S+|\nabla \psi|_\omega\sqrt{S}\leq \frac{3S}2+|\nabla \psi|^2_\omega
    \end{equation}
    where in the second inequality we applied Young's inequality. As for the second term, we have, by repeated application of Cauchy-Schwarz,
    \begin{equation*}
    |(II)|\leq \left(\sum_{i}a_i |v_i|^2\right)^{1/2}\cdot \left(\sum_{i}a_i\left|\sum_k \theta_{k\bar{i}}u_{\bar{k}}\right|^2\right)^{1/2}\leq C|\nabla v|_{\theta_u}\sqrt{ST}.
    \end{equation*}
    Applying Young's inequality for some $0<\epsilon\ll\min\{1,\delta\}$, we find 
    \begin{equation}\label{quant2}
        |(II)|\leq C|\nabla v|^2_{\theta_u}+\frac\epsilon4 ST.
    \end{equation}
    Hence, combining the bounds from \eqref{quant1} and \eqref{quant2}, we see
    \[
    \operatorname{Re}\left(\sum_{i,k} a_i v_iu_{k\bar{i}}u_{\bar{k}}\right)\leq \frac{3S}{2}+|\nabla \psi|_\omega^2+C|\nabla v|^2_{\theta_u}+\frac{\epsilon}4ST.
    \]
    We can estimate \eqref{estimate4.1} now as follows: 
    \begin{align*}
    0\geq\sum_{i,k}a_i|u_{k\bar{i}}|^2-CST-Ce^{Cj}S+S(h''|\nabla v|^2_{\theta_u}-np+\delta pT)-2p\operatorname{Re}\left(\sum_{i,k} a_i v_i u_{k\bar{i}}u_{\bar k}\right)
    \end{align*}
    which simplifies to 
    \begin{equation}\label{LQbound}
        0\geq (Sh''-C_B)|\nabla v|^2_{\theta_u}+S\left(\left(\frac{\delta p}{2}-C\right)T -(n+3)p\right)-C_Be^{Cj}\pare{1+S}.
    \end{equation}
    Constants denoted $C_B$ depend on $B$, while constants without the subscript do not. Choose $B,C$ sufficiently large so that both $B\kappa\ge L_1$ and  
    \[\frac{\delta p}{2}-C \ge \frac{\delta p}{4}\]
    hold. If $S\le C_B t^3$ or $S\le 4\abs{\nabla\psi}_\omega^2$, then $S\le Ce^{bj}$ since $1\le t\le C\pare{1+j}$ and we are done. After increasing the threshold constant, we may assume $Sh^{\prime\prime}-C_B\ge \frac{S}{t^3}$.

    Choose \[T_*\ge \max\brac{\frac{24C_Be^{Cj}}{B\delta},\frac{8(n+3)}{\delta}},\] growing exponentially in $j$. If $T\ge T_*$, then the positive trace term alone exceeds $C_B^{Cj}\pare{1+S}$, which is a contradiction. If $T<T_*$, every eigenvalue $\mu_i$ of $\theta_u$ relative to $\omega$ satisfies $\mu_i\le Ce^{Cj}T_*^{n-1}$, by Lemma \ref{determinantbounded}. Since $S>4\abs{\nabla \psi}_\omega^2$, we have \[\abs{\nabla v}_{\theta_u}^2\ge \frac{S}{4Ce^{Cj}T_*^{n-1}}.\] Dropping the nonnegative $T$ contribution and using $p\le B+1$ gives \[\frac{S^2}{4Ce^{Cj}T_*^{n-1}t^3}\le C_Be^{Cj}S.\] This makes $S$ exponentially bounded at the maximum in every case. The uniform upper bound for $h$ gives the desired bound for $Q$, and hence gives the estimate for $S$ claimed.
\end{proof}

\section{Boundary Hessian Estimates and Regularity}\label{section5}

In this section, we  move the analysis to the boundary of each $\Omega_j$. At a boundary point $p\in \partial \Omega_j$ and for $r>0$, denote the spherical cap $\Omega^{(r)}_j$ around $p$ by 
\[
\Omega^{(r)}_j = \Omega_j\cap B_r(p).
\]
We will again work in local holomorphic coordinates, considering the logarithmic version of the Monge-Amp\`ere equation
\begin{equation}\label{logdetmongeampere}
        \log \det\tilde{g} = \lambda_j u_j + F_j + \log\det g
\end{equation}
where $\tilde g$ denotes the metric associated to $\theta_u$. The main advantage of  \eqref{logdetmongeampere} is that differentiating yields the simpler second-order elliptic operator
\begin{equation}\label{linearmongeamp}
L=\tilde g^{i\bar{j}}\partial_i\partial_{\bar{j}}
\end{equation}
which is more amenable to both the geometry and the analysis of the estimates in this section. The methods in this section follow \cite{CKNS, GuanLi}.

\begin{remark}
    All quantities introduced up to this point that depend on $j$ satisfy a bound of the form $Q^{\pare{\ell}}_j\le A_\ell e^{B_\ell j}$ for some $A_\ell,B_\ell$ depending on the quantity $Q^{\pare{\ell}}$.
    
    We choose $\Lambda_j=C_\Lambda e^{c_\Lambda j}$ to dominate the data specified in Corollary \ref{wpolynomialbound}, Lemma \ref{backgroundgeo},  Lemma \ref{determinantbounded}, Proposition \ref{GradientEstimateuj} and $\|\nabla F_j\|_{L^\infty(X)}$. By $P(\Lambda_j)$, we mean a quantity bounded by $C\Lambda_j^m$ for some $C,m>0$, so every fixed polynomial expression in these bounds is bounded by $P(\Lambda_j)$.
\end{remark}

As in Lemma \ref{BoundaryC1Estimates}, we may apply a maximum principle-style argument to obtain a helpful comparison that we will use below.
\begin{lem}\label{approxboundary}
    Let $p\in \partial \Omega_j$, and let $r>0$ satisfy  $0<\frac1r\leq P(\Lambda)$. Let $G\geq 0$ be a function, and let $C'$ satisfy $0<C'\leq P(\Lambda)$ and 
    \[
    L|z-p|^2\le C'G.
    \]
    Suppose that $V\in C^2(\Omega^{(r)}_j)$ is a nonnegative barrier function such that $V(p)=0$ and 
    \begin{equation}
        L(V)\leq -G, \qquad V_{\nu(p)}(p)\leq P(\Lambda)
    \end{equation}
    where $\nu(p)$ denotes the inward normal direction at $p$. If a real-valued function $W\in C^2(\Omega_j^{(r)})\cap C^1(\overline{\Omega_j^{(r)}})$ satisfies
    \begin{align*}
        L(W)\ge -H_1G,\qquad \abs{W}\le H_2
    \end{align*}
    and on $\del\Omega_j\cap B_r(p)$
    \begin{align*}
        W(p)=0,\qquad W\le H_3\abs{z-p}^2
    \end{align*}
    for $0\leq H_1,H_2,H_3\le P\pare{\Lambda}$, then 
    \[W_{\nu(p)}(p)\le P(\Lambda).\]
\end{lem}
\begin{proof}
    Choose $B= H_3+H_2r^{-2}\leq P(\Lambda)$. On the actual boundary, $B\abs{z-p}^2-W\ge 0$ and on the artificial sphere boundary the same bound holds because $Br^2\ge H_2$. Since $V\ge 0$, $AV+B\abs{z-p}^2-W\ge 0$ on $\del\Omega_r$ for all $A\ge 0$. Define
    \begin{align*}
        A= H_1+C'B\leq P(\Lambda).
    \end{align*}
    We thus have 
    \begin{align*}
        L\pare{AV+B\abs{z-p}^2-W}\le 0.
    \end{align*}
    The maximum principle implies that $W\le AV+B\abs{z-p}^2$ on all of $\Omega_j^{(r)}$. Both sides agree at $p$, so differentiation in the inward direction yields 
    \begin{equation*}
        W_{\nu(p)}(p)\le AV_{\nu(p)}\pare{p}\le P(\Lambda).\qedhere
    \end{equation*}
\end{proof}

\begin{prop}\label{HessianBoundaryControl}
    For sufficiently large $j$, there are uniform constants $C,L_2>0$ such that $\sup_{\partial\Omega_j}\abs{\nabla^2 u_j}_\omega\le Ce^{L_2j}$.
\end{prop}

\begin{proof}
As before, to simplify the notation, we suppress the subscripts $j$. Recall that $\theta_w\ge \tau\omega$ with $\tau^{-1}\le P\pare{\Lambda}$. Setting
\[s=u-w,\] 
we see that $\norm{s}_{C^1}\le P\pare{\Lambda}$. At any boundary point $p\in \partial \Omega$, let $x$ denote the inward real normal direction and $\xi,\xi^\prime$ denote real tangential directions. To estimate the real Hessian of $u$, we compute the tangential, the mixed, and the normal components separately.

\textit{Tangential Hessian:} Since $s$ vanishes on $\partial \Omega$, if $\rho$ is a defining function for $\Omega$ then there exists a function $H$ such that 
\[s=H\rho.\]
Let $p\in \partial \Omega$, and choose coordinates such that $\rho_\xi(p)=\rho_{\xi'}(p)=0$ and $\rho_{x}(p)=-1$. Differentiating in the inward normal direction $x$, we see
\[
s_{x}(p)=H_{x}(p)\rho(p)+H(p)\rho_{x}(p) = -H(p).
\]
Differentiating in the real tangential directions $\xi, \xi'$ at $p$, we obtain the following identity in these coordinates:
\begin{align*}
    s_{\xi\xi^\prime}(p) &= H(p)\rho_{\xi\xi'}(p) + H_\xi(p)\rho_{\xi'}(p)+H_{\xi'}(p)\rho_{\xi}(p)+H_{\xi\xi'}(p)\rho(p)\\
    &=H(p)\rho_{\xi\xi^\prime}(p) \\
    &=-s_{x}(p)\rho_{\xi\xi^\prime}(p).
\end{align*}
By the gradient estimate for $s$ from Proposition \ref{GradientEstimateuj} and Corollary \ref{wpolynomialbound} along with the boundary geometry from Lemma \ref{backgroundgeo}, we have
\[
\abs{s_{\xi\xi^\prime}(p)}\le P\pare{\Lambda}.
\]
In particular, Corollary \ref{wpolynomialbound} controls the derivatives of $w$, so we may also conclude
\begin{equation}\label{realtanghessian}
    \abs{u_{\xi\xi^\prime}(p)}\le P\pare{\Lambda}.
\end{equation}
We also record that since $s\geq 0$ by Lemma \ref{ujclose}, these estimates also imply that
\begin{equation}\label{diffxbound}
0\leq s_x(p)\leq C\Lambda  .  
\end{equation}

As stated at the start of this section, let $\tilde{g}_{i\bar{j}} = \theta_{i\bar{j}} +u_{i\bar{j}}$ and consider the operator $L$ from  \eqref{linearmongeamp}. Let
\begin{align*}
    T &= \operatorname{Tr}_{\theta_u}\omega.
\end{align*}
If $d$ is the inward distance to the boundary, in the spherical cap $\Omega^{(r)}$, define the barrier
\begin{align*}
    V(d) &= s+td-Nd^2.
\end{align*}
Direct calculation gives
\begin{align*}
    LV &\le n-\tau T+C\Lambda\pare{t+CNr}T-2N\abs{\del d}^2_{\theta_u}.
\end{align*}
Choose $t=\frac{c\tau}{\Lambda}$ and $N=C\pare{1+\Lambda}\tau^{-n}$, with $c$ small, $C$ large, and $r$ sufficiently small satisfying $\frac1r\leq P(\Lambda)$ and so that 
\begin{align*}
    r\le \ell, \qquad Nr\le \frac{t}{2}, \qquad C\Lambda\pare{t+2Nr}\le \frac{\tau}{2}.
\end{align*}
Applying the arithmetic-geometric mean inequality to $\tilde{g}^{-1}\pare{\frac{\tau}{4}g+N\del d\otimes \bar{\del} d}$, we have
\begin{align*}
    \frac{\tau}{4}T +N\abs{\del d}^2_{\theta_u} &\ge c\pare{\frac{\tau^{n-1}N}{\det\pare{\omega^{-1}\theta_u}}}^{1/n}.
\end{align*}
Combining this inequality with our bound on $LV$, we obtain
\begin{equation*}
    LV \le -c_0\tau\pare{1+T}
\end{equation*} 
for some constant $c_0>0$.
Furthermore, $V\ge 0$ on $\Omega_r$, $V(p)=0$, and \[V_x(p)=s_x(p)+t\le P\pare{\Lambda}.\]

\textit{Mixed Hessian:} For the mixed derivatives, we consider the logarithmic form in  \eqref{logdetmongeampere}.
 For a fixed real coordinate direction $\xi$, differentiation gives 
    \begin{equation*}
        L\pare{\del_\xi s} = \lambda \del_\xi u + \del_\xi F + \del_\xi \log\det g - \tilde{g}^{i\bar{j}}\del_\xi(\theta_w)_{i\bar{j}}.
    \end{equation*}
    Using the bounds above, we see
    \begin{align}\label{xisbound}
        \abs{L\pare{\del_\xi s}} &\le P\pare{\Lambda}\pare{1+T}.
    \end{align}
    Now consider the tangential vector field
    \begin{align*}
        G_\xi = \del_\xi-\frac{\rho_\xi}{\rho_x} \del_x,
    \end{align*}
    defined so that $G_\xi(\rho)=0$. We see that $G_\xi(s)=0$ on $\del\Omega$.
    Setting $z_n=x+\sqrt{-1}y$, we see that 
    \begin{equation}\label{directionidentity}
    s_{x\bar{j}} = 2s_{n\bar{j}}+\sqrt{-1}s_{y\bar{j}}.    
    \end{equation}
    Applying $L$ to $G_\xi (s)$, we get  
    \begin{align*}
        \abs{L\pare{G_\xi (s)}} &= \abs{L(s_\xi) -\frac{\rho_\xi}{\rho_x} L(s_x) -\pare{L\pare{\frac{\rho_\xi}{\rho_x}}}s_x -2\operatorname{Re}\pare{\tilde{g}^{i\bar{j}}\pare{\frac{\rho_\xi}{\rho_x}}_i s_{x\bar{j}}}}.
        \end{align*}
    Applying the triangle inequality along with the decomposition in \eqref{directionidentity} followed by Young's inequality, we have     
    \begin{align*}\abs{L\pare{G_\xi (s)}} &\le P\pare{\Lambda}\pare{1+T} + \abs{2\operatorname{Re}\pare{\tilde{g}^{i\bar{j}}\pare{\frac{\rho_\xi}{\rho_x}}_i (\sqrt{-1}s_{y\bar{j}})}} \\
    &\le P\pare{\Lambda}\pare{1+T}+\abs{\nabla \left(\frac{\rho_\xi}{\rho_x}\right)}^2_{\theta_u} + \abs{\nabla s_y}^2_{\theta_u}.
    \end{align*}
    Finally, recalling how we defined $T$, we may absorb the second term in the previous expression to arrive at
    \begin{align*}
    |L(G_\xi(s))|&\le P\pare{\Lambda}\pare{1+T} + \abs{\nabla s_y}^2_{\theta_u}.
    \end{align*}
    On the other hand, we see
    \begin{align*}
        L\pare{s_y^2} &= 2\abs{\nabla s_y}^2_{\theta_u}+2s_yL(s_y)\ge 2\abs{\nabla s_y}^2_{\theta_u} -P\pare{\Lambda}\pare{1+T}.
    \end{align*}
    Thus $L\pare{s_y^2\pm G_\xi(s)}\ge -P\pare{\Lambda}\pare{1+T}$. On the actual boundary, we have 
    \[s_y=\frac{\rho_y}{\rho_x}s_x.\]
    Given our choice of defining function $\rho$, we see that $s_y\pare{p}=0$. Furthermore, since
    \[s_y^2\le P\pare{\Lambda}\abs{z}^2,\]
    and we have polynomial bounds on $s_y^2\pm G_\xi (s)$ by Corollary \ref{wpolynomialbound}, we may apply Lemma \ref{approxboundary} to get
    \begin{align*}
        \abs{\del_x\pare{G_\xi (s)}\pare{p}}=\abs{s_{\xi x}\pare{p}-\pare{\frac{\rho_\xi}{\rho_x}}_x(p)s_x\pare{p}} \le P\pare{\Lambda}.
    \end{align*}
    In particular, we can rewrite this bound as follows:
    \begin{equation}\label{mixedhessian}
        \abs{u_{\xi x}\pare{p}} \le P\pare{\Lambda}.
    \end{equation}

    We would like to recover the normal component from the determinant of $\tilde{g}$, and thus need to argue that the complex tangential block does not degenerate (if $n=1$, this argument is unnecessary, so we assume that $n\ge 2$). For $p\in \partial \Omega$, define the set 
    \[
    U_p=\{\zeta\in T^{1,0}_p(\partial \Omega): |\zeta|_\omega=1\}
    \]
    and consider 
    \begin{align*}
        m_0 &= \min_{p\in\del\Omega}\pare{\min_{\zeta\in U_p} \pare{\theta_u\pare{\zeta,\bar{\zeta}}}}.
    \end{align*}
    If $m_0\ge \frac{\tau}{2}$, then our lower bound already holds. Otherwise, consider a minimizing point $p$ and direction, and normalize coordinates so that this direction is the unit normalization of $\del_{z_1}$. Consider the unit tangential vector field on $\Omega^{(r)}$ given by 
    \[
    \zeta = \frac{\rho_n\partial_{z_1}-\rho_1\partial_{z_n}}{|\rho_n\partial_{z_1}-\rho_1\partial_{z_n}|_\omega}.
    \] 
    Define the barrier $\Phi$ on $\Omega^{(r)}$ by
    \[
    \Phi = \pare{\theta_w}_{i\bar{j}}\zeta^i\bar{\zeta}^j+\frac{\rho_{i\bar{j}}}{\rho_x}\zeta^i\bar{\zeta}^js_x -m_0.
    \]
    On $\partial \Omega\cap B_r(p)$, since 
    \[
    \theta_u(\zeta,\bar{\zeta}) = \pare{\theta_w}_{i\bar{j}}\zeta^i\bar{\zeta}^j+\frac{\rho_{i\bar{j}}}{\rho_x}\zeta^i\bar{\zeta}^js_x
    \]
    we consequently have $\Phi\geq 0$ on $\partial \Omega\cap B_r(p)$. At the minimizing point $p_0$ in the definition of $m_0$, we have that $\Phi(p_0)=0$ and
    \begin{equation}\label{b0bound}
    -\frac{\rho_{i\bar{j}}}{\rho_x}\zeta^i\bar{\zeta}^js_x = \pare{\theta_w}_{i\bar{j}}\zeta^i\bar{\zeta}^j -m_0\ge \frac{\tau}{2}.   
    \end{equation}
    By \eqref{diffxbound}, we hence have 
    \[-\frac{\rho_{i\bar{j}}}{\rho_x}\zeta^i\bar{\zeta}^j\ge \frac{c\tau}{\Lambda}.\] 
    Using the polynomial bounds and repeating the mixed derivative argument, we obtain
    \begin{align*}
        L(\Phi) &\le \abs{\nabla s_y}^2_{\theta_u} + P\pare{\Lambda}\pare{1+T}.
    \end{align*}
    Thus 
    \[
    L\pare{s_y^2-\Phi}\ge -P\pare{\Lambda}\pare{1+T}.\] 
    On the true boundary, $\Phi\ge 0$ and $s_y^2\le P\pare{\Lambda}\abs{z-p_0}^2$, while on the spherical boundary, we have $\abs{s_y^2-\Phi}\le P\pare{\Lambda}$ by the coefficient and gradient bounds.
    Again, applying Lemma \ref{approxboundary}, we see 
    \[
    \Phi_x(p_0)\geq -P(\Lambda)
    \]
    meanwhile differentiation yields
    \[\Phi_x(p_0)=\frac{\rho_{i\bar j}(p_0)}{\rho_x(p_0)}\zeta^i\zeta^{\bar{j}}s_{xx}(p_0)+O(P(\Lambda)).\]
    The inequality \eqref{b0bound} in turn implies that $s_{xx}(p_0)\leq P(\Lambda)$, and hence $u_{xx}\le P\pare{\Lambda}$.

    \textit{Normal Components:} In the branch $m_0<\tau/2$, the preceding normal estimate bounds the eigenvalues at $p_0$ above, so the determinant lower bound gives $m_0\ge c\Lambda^{-N_0}$. In the branch $m_0\ge\tau/2$, the same bound follows directly from the polynomial control of $\tau^{-1}$. Thus in either branch the complex tangential block has this lower bound at every boundary point.
    Now we can prove the bound for the normal direction. At an arbitrary point $q\in \partial \Omega$, use an $\omega$-unitary frame whose first $n-1$ vectors are complex tangential and decompose $\tilde{g}$ as
    \begin{align*}
        \pare{\tilde{g}_{i\bar{j}}} &= \begin{pmatrix}
            A & \beta \\ \beta^* & c
        \end{pmatrix}
    \end{align*}
    where $A$ is an $(n-1)\times (n-1)$ matrix, $\beta$ is an $(n-1)\times 1$ vector, and $c$ is a scalar entry. 
    We see that $A$ obeys the bound 
    \[A\ge \tilde{c}_0\Lambda^{-N_0}I_{n-1},\] 
    and our mixed and tangential estimates give polynomial bounds for $\beta$. Therefore, the Schur-complement identity gives
    \begin{align*}
        c &=\frac{\det \pare{\tilde{g}_{ij}}}{\det A} + \beta^*A^{-1}\beta \le P\pare{\Lambda}.
    \end{align*}
    Again because $u_{n\bar{n}}=\frac{1}{4}\pare{u_{xx}+u_{yy}}$, the positivity and tangential bounds for $u_{yy}$ give both bounds for $u_{xx}$. Therefore,
    \begin{align*}
        \sup_{\del\Omega}\abs{\nabla^2 u}_\omega &\le P\pare{\Lambda}.
    \end{align*}
    The conversion from coordinate to covariant derivatives adds only controlled first-derivative terms. 
    
    Reintroducing the exhaustion index, since $\Lambda_j\le Ce^{Cj}$, we indeed have 
    \begin{align*}
        \sup_{\del\Omega_j}\abs{\nabla^2 u_j}_\omega &\le Ce^{L_2j}.\qedhere
    \end{align*}
\end{proof}

The boundary Hessian bound may grow with the exhaustion level. At a boundary maximum of the trace test function, however, $u_j-\psi\ge \kappa j-C$ offsets that growth. The interior argument distinguishes the two available controls for the log density.

\begin{lem}\label{TraceUpperBound}
    For sufficiently large $j$, there exist uniform constants $B_1,C>0$ such that
    \begin{itemize}
        \item[(i)] under the hypotheses of Theorem \ref{maintheorem} (i), $\operatorname{Tr}_\omega\theta_{u_j}\le Ce^{-B_1\psi}$ on $\Omega_j$;
        \item[(ii)] under the hypotheses of Theorem \ref{maintheorem} (ii), $\operatorname{Tr}_\omega\theta_{u_j}\le Ce^{-F_j^- -B_1\psi}$ on $\Omega_j$.
    \end{itemize}
\end{lem}

\begin{proof}
    To simplify the notation, we suppress the subscripts $j$ and set 
    \[
    H = \operatorname{Tr}_\omega \theta_u,\qquad v=u-\psi,\qquad m=\inf_\Omega v.
    \]
    For item (i), the goal is to maximize the quantity 
    \[Q=\begin{cases}
        \log H-Av &\text{in case (i)} \\
        \log H-Av+F_j^- & \text{in case (ii)}
    \end{cases}\] 
    Fix \(A>0\) sufficiently large to satisfy all the requirements below, and let $x_0\in\overline\Omega$ maximize $Q$. If $x_0\in \partial \Omega$ then by Proposition \ref{HessianBoundaryControl}, $\log H\le L_2 j + C$. Since $\psi(z)=-j$ for all $z\in \partial \Omega$, we have that $v\ge \kappa j-C$. Additionally, in the case of (ii), we have a uniform upper bound on $F_-$. Hence, in either case 
    \begin{align*}
        Q\pare{x_0} \le C_A+\pare{L_2-A\kappa}j
    \end{align*}
    provided $A\kappa\ge L_2$. 
    
    If instead $x_0$ is an interior maximum and $H\pare{x_0}<Ce^{-L_3\psi\pare{x_0}}$ (where $L_3\ge B_2$ as in Lemma \ref{Fapprox} (i) in the first setting or $L_3=L_2$ in the second setting), then $Q\pare{x_0}\le \pare{A\kappa-L_3}\psi\pare{x_0}+C_A\le C_A$, provided $A\kappa \ge L_3$. Therefore, we can assume that $x_0$ is an interior maximizer satisfying $H\pare{x_0}\ge Ce^{-L_3\psi\pare{x_0}}$.

    Choose normal holomorphic coordinates for $\omega$ at the point $x_0$ which diagonalize $\theta_u$ as in Proposition \ref{GradientEstimateuj}. Computing $LQ$ directly, we find
    \begin{equation}\label{Lbound}
        L\pare{\log H - Av} \ge \frac{\Delta_\omega F}{H}-CT-A\pare{n-\delta T}.
    \end{equation}
    Meanwhile, computing $LH$ and using the fact that $\theta$ is closed gives
    \begin{align*}
        LH \ge \sum_{i,k}a_i\del_k\del_{\bar{k}}\tilde{g}_{i\bar{i}}-CHT.
    \end{align*}
    By considering \eqref{logdetmongeampere}, we may differentiate twice to refine our lower bound to 
    \begin{align*}
        LH &\ge \Delta_\omega\pare{\lambda u +F}+\sum_{i,k,\ell}a_ia_\ell\abs{\del_k \tilde{g}_{i\bar{\ell}}}^2-CHT-C.
    \end{align*}
    Weighted Cauchy-Schwarz on $\del_i H= \sum_k\del_k\tilde{g}_{i\bar{k}}$ yields
    \begin{align*}
        L\pare{\log H - Av} &\ge  \frac{LH}{H}-\frac{\abs{\nabla H}^2_{\theta_u}}{H^2} -A\pare{n-\delta T} \\
        &\ge \frac{\Delta_\omega\pare{\lambda u +F}}{H}-CT-\frac{C}{H} -A\pare{n-\delta T}.
    \end{align*}
    Noting that $\Delta_\omega u=H-\operatorname{Tr}_\omega \theta$, $HT\ge n^2$, and that $0\le\lambda\le 1$, we can rewrite the bound as \eqref{Lbound}. We now split into our two cases:
    \begin{itemize}
        \item[(i)] Equation \eqref{Lbound} implies
        \[
        0\ge LQ\pare{x_0}\ge \pare{A\delta-C}T\pare{x_0}-An-C.
        \]
        Choosing $A$ so that $A\delta>C$ we find that $T\pare{x_0}\le C_A$. By Lemma \ref{Fapprox} (i)
        \[
        H\pare{x_0} \le n\det(\omega^{-1}\theta_u\pare{x_0})T^{n-1}\pare{x_0} \le C_A e^{-B\psi\pare{x_0}}.
        \]
        By our choice of $A$, we have $A\kappa\ge B$, so we see that $Q\pare{x_0}\le C_A$, and thus $Q\le C_A$ on all of $\Omega$. Using $v=u-\psi$, we conclude $H\le Ce^{-A\psi}$.

        \item[(ii)] Observe that 
        \[
        \frac{\operatorname{Tr}_\omega \pare{\pot F_j^-+C_*\omega}}{H} \le \operatorname{Tr}_{\theta_u}\pare{\pot F_j^-+C_*\omega}.
        \]
        Expanding and dropping the nonpositive term $-\frac{nC_*}{H}$, we find
        \[
        \frac{\Delta_\omega F_j^-}{H} \le LF_j^-+C_*T.
        \]
        Since $\Delta_\omega F^+_j\geq -nC_*$ and $HT\geq n^2$, we obtain 
        \[
        L\pare{\log H + F_j^-}\ge -C^\prime T.
        \]
        Returning to $Q$, this shows
        \[
        0\ge LQ\pare{x_0}\ge \pare{A\delta -C^\prime}T\pare{x_0}-An.
        \]
        Again, for $A$ sufficiently large, $T\pare{x_0}\le C_A$. Therefore,
        \begin{align*}
        H\pare{x_0}e^{F_j^-\pare{x_0}}&\le n\det\omega^{-1}\theta_{u}\pare{x_0} e^{F_j^-\pare{x_0}}T\pare{x_0}^{n-1}\\
        &\le ne^{\lambda u\pare{x_0}+F_j^+\pare{x_0}}T\pare{x_0}^{n-1}\\
        &\le C_A.
        \end{align*}
        This similarly bounds $Q\pare{x_0}$, and thus $H\le Ce^{-F_j^- -A\psi}$.\qedhere
    \end{itemize}
\end{proof}

On a fixed compact subset of $U$, the barrier weight and, in the case of (ii), $e^{-F_j^-}$ are uniformly bounded. Thus, the trace estimate and the positive local density bound give the ellipticity needed for the interior H\"older estimates.

\begin{prop}\label{UniformLocalRegularity}
    Let $0<\alpha<\beta$. Then for each compact pair $K_1\Subset K_2\Subset U$, there is a constant $C_{K_1,K_2,\alpha}$ such that for all sufficiently large $j$, $\norm{u_j}_{C^{2,\alpha}\pare{K_1}}\le C_{K_1,K_2,\alpha}$.
\end{prop}
\begin{proof}
    Cover $K_1$ by finitely many fixed holomorphic coordinate balls $B_\nu\Subset B_\nu^\prime\Subset B_\nu^{\prime\prime}\Subset\operatorname{int}K_2$. These balls lie in $\Omega_j$ for all sufficiently large $j$. On $B_\nu^{\prime\prime}$, choose a smooth local potential $a$ for $\theta$ and put $\hat{u}_j=u_j+a$. By Lemma \ref{ujclose} and Lemma \ref{TraceUpperBound}, $\widehat u_j$ is uniformly bounded there and its Euclidean Laplacian is nonnegative and uniformly bounded. Interior Poisson estimates give a uniform gradient bound on $B_\nu^\prime$. On that same ball, the density $e^{\lambda_ju_j+F_j}\det(g_{k\bar\ell})$ is uniformly bounded below by a positive constant, and its $n$-th root has a uniform $C^{0,\alpha}$ norm by Lemma \ref{Fapprox} and the gradient bound. Apply Wang’s interior estimate \cite[Theorem 1.2]{Wang} to $\hat{u}_j$ on $B_\nu'$ to obtain a uniform $C^{2,\alpha}$ estimate on $B_\nu$. Subtracting $a$ and taking the finite cover proves the result.
\end{proof}

\begin{cor}\label{locallysmoothu}
    If $F$ is smooth, then for any  $K_1\Subset K_2\Subset U$, $0<\alpha<1$, $k\ge 2$, and sufficiently large $j$, there is a uniform constant $C_{K_1,K_2,k,\alpha}$ such that $\norm{u_j}_{C^{k,\alpha}\pare{K_1}}\le C_{K_1,K_2,k,\alpha}$.
\end{cor}
\begin{proof}
    Let $a$ be a local potential for $\theta$ and set $\hat{u}_j=u_j+a$. Differentiating 
    \begin{align*}
        \log\det \pare{\hat{u}_j}_{i\bar{k}} &= \lambda_j\hat{u}_j+F_j-\lambda_j a +\log\det g_{i\bar{k}}
    \end{align*}
    in a real coordinate direction $\mu$ and setting $A_j^{i\bar{k}}=\pare{\pare{\hat{u}_j}_{i\bar{k}}}^{-1}$ yields the linear equation
    \begin{align*}
        A_j^{i\bar{k}}\del_{i\bar{k}}\pare{\del_\mu \hat{u}_j}-\lambda_j\del_\mu\hat{u}_j &= \del_\mu F_j-\lambda_j\del_\mu a +\del_\mu\log\det g_{i\bar{k}}.
    \end{align*}
    Note that $\abs{\lambda_j}\le 1$ and the operators are uniformly elliptic with uniformly $C^\alpha$ principal coefficients. By Corollary \ref{smoothFapprox}, interior Schauder estimates then give a $C^{3,\alpha}$ bound on a smaller neighborhood, and with iteration we get all higher bounds. A finite covering argument proves the result for all compact pairs contained in $U$.
\end{proof}

With this, we are now able to prove Theorem \ref{maintheorem}.
\begin{proof}[Proof of Theorem \ref{maintheorem}]
    We know that $\varphi_j\to\varphi_\infty$ in $L^1\pare{X}$ and locally uniformly by Remark \ref{limitworks}. Further, Lemma \ref{ujclose} asserts that $u_j\to \varphi_\infty\big|_U$ locally uniformly in $U$.

    Now fix $0<\alpha<\alpha^\prime<\beta$. By Proposition \ref{UniformLocalRegularity}, we have uniform $C^{2,\alpha^\prime}$ bounds on fixed compact subsets. Compact embedding on smaller coordinate balls and diagonal extraction give local $C^{2,\alpha}$ convergence of $u_j\to \varphi_\infty\big|_U$. Thus, $\varphi_\infty\in C^{2,\alpha}_{\loc}\pare{U}$ for all $0<\alpha<\beta$. Equivalently, for each such $\alpha$, the chosen subsequence converges locally in $C^{2,\alpha}$ without another extraction by precompactness and uniqueness of its local uniform limit.

    On each fixed compact subset of $U$, $\lambda_j u_j\to 0$ and $F_j\to F$ uniformly. Passing to the limit gives the classical equation on $U$. Remark \ref{limitworks} already established the global non-pluripolar equation for $\varphi_\infty$ on $X$. Subtracting a constant to achieve the normalized solution affects neither the regularity nor the defining equation.

    If $F$ is smooth on $U$, the limit equation bootstraps locally. Its complex Hessian is positive definite, as it is semipositive and has positive determinant. Consequently, it is uniformly elliptic on each smaller compact coordinate ball. Differentiating the limit equation and iterating the interior Schauder estimates gives smoothness on all of $U$. This last step can be carried out directly on the limit equation.
\end{proof}

\begin{proof}[Proof of Corollary \ref{Answer}]
    In the setting of Theorem \ref{maintheorem}, we take $U=\amp$ and $Z$ to be the non-K\"ahler locus of $\sqbr{\theta}$. By \cite[Section 1.4]{BEGZ}, the class $\sqbr{\theta}$ contains a K\"ahler current $\theta+dd^c\phi$ with analytic singularities precisely along $Z$. On a log resolution of its defining ideal, chosen to be an isomorphism over $U$, the pullback of $\phi$ is locally a sum of logarithms of simple normal-crossing divisor coordinates plus a smooth remainder. After subtracting a constant, we may assume $\phi\le -1$. Thus $\phi$ is a barrier as in Theorem~\ref{maintheorem}. Since the density is smooth and strictly positive on $X$, its logarithm has bounded Laplacian, so case~(i) gives smoothness on $U$.
\end{proof}

\begin{remark}\label{Fujiki}
To see that Corollary \ref{Answer} also extends to the case when $X$ belongs to Fujiki class $\mathcal{C}$, let $(Y,\hat\omega)$ be a compact K\"ahler manifold. Consider the modification
\[
\pi:Y\to X
\]
which is an isomorphism over $\amp$, obtained by resolving a K\"ahler current $T_0\in [\theta]$ with analytic singularities precisely outside $\amp$. In particular, $\pi^*T_0=\beta+[D]$, where $\beta$ is smooth and $D$ is a simple normal-crossing divisor whose support is $Y\setminus \pi^{-1}(\amp)$. If $\omega_X$ is a Hermitian metric on $X$ and $T_0\ge c\omega_X$, then $\beta\ge c\pi^*\omega_X$ off $D$, and thus everywhere by smoothness. By \cite{DemaillyPaun}, consider an effective exceptional divisor $E'$ and a smooth representative $\eta\in [E']$ such that $\hat\omega=\beta-\eta$ is K\"ahler. Thus, $\hat\omega+[D+E']=\pi^*\theta+\pot \hat\phi$, where $\hat\phi$ is our barrier function on $Y$ possessing analytic logarithmic normal-crossing singularities precisely on $Y\setminus \pi^{-1}(\amp)$, which we normalize by $\hat{\phi}\le -1$. In particular, $\pi^{-1}\pare{\amp}\subseteq \operatorname{Amp}\pare{\pi^*\theta}$.

Considering the pullback $\pi^*(f\omega^n) = g\hat{\omega}^n$, we see that $g$ vanishes along each divisorial component of the ramification locus of $\pi$. However, over $\pi^{-1}(\amp)$, $g$ is smooth and positive. Locally, we may represent $g$ as
\[
\log g=\log|\det d\pi|^2+h
\]
for some function $h:Y\to \mathbb{R}$ which is smooth in a coordinate neighborhood across $E$. Thus
\[
\Delta_{\hat\omega} \log(g(z))\geq -C
\]
distributionally for some constant $C>0$. Hence, we are in condition (i) of Theorem \ref{maintheorem}.

Let $\varphi$ denote the normalized solution downstairs. By compatibility of the non-pluripolar product with modifications,
\[
\inner{\pare{\pi^*\theta+\pot\pare{\varphi\circ\pi}}^n}=\pi^*\pare{f\omega^n}=g\hat\omega^n.
\]
We see
\[
\int_Y g\hat{\omega}^n = \int_X f\omega^n = \vol = \operatorname{Vol}\pare{\pi^*\theta}.
\]
Thus, $\varphi\circ\pi$ is the full-mass equation in $\pi^*[\theta]$ up to a constant. Applying Theorem \ref{maintheorem} to $\varphi\circ\pi$ gives smoothness on $\pi^{-1}(\amp)$. Since $\pi$ is biholomorphic on $\pi^{-1}(\amp)$, smoothness descends to $\varphi$ on $\amp$. 
\end{remark}

\appendix

\section{Heat Kernel Estimates and Regularization}\label{heatappendix}

We will use $P_t=e^{t\Delta_\omega}$ to denote the scalar heat operator and $Q_t$ to be the Hodge heat operator on $(1,1)$-forms, with the same time normalization as on functions. By the K\"ahler identities, the Hodge Laplacian preserves type and commutes with $\del$ and $\bar{\del}$. Consequently, the identity
\[\pot \pare{P_t u}=Q_t\pare{\pot u}\]
holds for smooth functions and, by duality, for distributions. Because $\omega$ is harmonic and $X$ is K\"ahler, $Q_t\omega=\omega$. We use $K_t\pare{x,y}$ for the scalar kernel and $\mathcal{K}_t\pare{x,y}$ for the Hodge kernel. We take $0<t_0\le \min\brac{1,\operatorname{inj}\pare{\omega}^n}$ so that all operations below are valid.

\begin{lem}\label{heatbounded}
    For $0<t<t_0$, we have $0\le K_t\pare{x,y}$ and for $k=0,1,\ldots$,
    \[
    \abs{\nabla_x^k K_t\pare{x,y}}\le C_kt^{-n-k/2}e^{-c_kd\pare{x,y}^2/t}.
    \]
    The same bound holds for derivatives in $y$.
\end{lem}
\begin{proof}
    The Gaussian bound on $K_t$ is due to the general heat-kernel theory of Aronson and Grigor'yan \cite{Aronson, Grigor'yan1}. For the first spatial derivative, we use the gradient estimate for $\log K_t$ in \cite{Sheu}, and multiply by the bound on $K_t$. The higher spatial derivatives follow from the estimates of \cite{Hsu}. The corresponding bounds for derivatives in $y$ follow from self-adjointness.
\end{proof}

\begin{lem}\label{chainrule}
    Suppose $a$ is continuous and $\Delta_\omega a\ge -q$ distributionally on an open set, where $q\ge 0$ is continuous. Then for smooth, increasing, convex functions $h$, 
    \[
    \Delta_\omega h(a)\ge -h'(a)q.
    \]
\end{lem}
\begin{proof}
    On a compact subdomain, let $\zeta\ge 0$ be a cutoff function equal to one on a larger neighborhood and set $a_t=P_t\pare{\zeta a}$. In distributions, 
    \[\Delta_\omega\pare{\zeta a}\ge -\zeta q+\sqbr{\Delta,\zeta}a.\]
    The commutator distribution has order at most one and is supported away from the compact subdomain. After pairing with the heat kernel, the off-diagonal bounds for the kernel and its first derivative show its contribution is $o(1)$ uniformly there. Thus, $a_t\to a$ uniformly and 
    \[\Delta_\omega a_t\ge -q-o(1).\] 
    Returning to $h$, we see 
    \[\Delta_\omega h\pare{a_t}\ge -h^\prime\pare{a_t}\pare{q+o(1)}\]
    on the compact subdomain, and we can pass to the distributional limit using the uniform convergence of $a_t$ and $h^\prime\pare{a_t}$.
\end{proof}

\begin{lem}\label{heathessian}
    If $\pot u\ge -C_0\omega$ on $X$ in the sense of currents, then there is a uniform constant $C_1$ such that $\pot P_t u\ge -C_1\omega$ for $0<t<t_0\le 1$.
\end{lem}
\begin{proof}
    Let 
    \[S=\pot u+C_0\omega\]
    be a closed, positive current. Its trace measure $\nu$ is the measure $nS\wedge \omega^{n-1}$. In a local unitary frame, the coefficients of $S$ are measures dominated in total variation by a fixed multiple of $\nu$. Equivalently $S=A\, d\nu$, where $A$ is a measurable positive semidefinite hermitian form of trace one, $\nu$-a.e. This follows from the Radon-Nikodym theorem applied to the coefficient matrix, using positivity to dominate the off-diagonal coefficients by the diagonal ones.

    For sufficiently small $r_0>0$, there is a constant $C_S$ such that
    \[
    \nu\pare{B_r\pare{x}}\le C_Sr^{2n-2}
    \]
    for all $x\in X$ and $0<r\le r_0$. To obtain a uniform bound in $x$, assume $n\ge 2$ and choose a finite cover by smaller coordinate charts contained inside larger charts. Coordinate distance and volume-form traces are uniformly comparable on the relevant larger charts. In a chart, apply the Jensen-Lelong mass monotonicity formula to $S$ on Euclidean balls centered at the coordinate of $x$ (the integral of $S\wedge\pare{\pot \abs{z}^2}^{n-1}$ divided by $r^{2n-2}$ is nondecreasing with radius). Compare with a fixed larger radius available for every center in the smaller chart, and bound the larger-ball mass by a constant times $\nu\pare{X}$. The $n=1$ case follows from the finiteness of $\nu$.

    By the short-time parametrix for $\mathcal{K}_t\pare{x,y}$ \cite[Theorem 2.30]{BGV}, we can write $\mathcal{K}_t=L_t+E_t$, where $L_t$ is a nonnegative scalar multiple of the Levi-Civita parallel transport near the diagonal, cut off to zero outside that neighborhood, and $E_t$ is an error term. Since $\omega$ is K\"ahler, parallel transport preserves positive $(1,1)$-forms. Moreover, taking a sufficiently long expansion gives
    \[
    \norm{E_t\pare{x,y}} \le Ct^{1-n}e^{-cd\pare{x,y}^2/t}+Ct\, \qquad 0<t\le t_0.
    \]

    We see that $\int L_t S\ge0$, so we only need to bound the integral of the error term. We split the integral over $B_{\sqrt{t}}\pare{x}$, the annuli $2^k\sqrt{t}\le d\pare{x,y}<2^{k+1}\sqrt{t}$ for $k=0,1,\ldots$ remaining in $B_{r_0/2}\pare{x}$, and the portion outside this ball. Under this decomposition, we find
    \[
    \int_X e^{-cd\pare{x,y}^2/t}\, d\nu\pare{y}\le C_S t^{n-1}\pare{1+\sum_{k=0}^\infty 2^{\pare{k+1}\pare{2n-2}}e^{-c4^k}}+e^{-c^\prime r_0^2/t}\nu\pare{X}\le C_S^\prime t^{n-1}.
    \]
    Thus, 
    \[\norm{\int E_t\pare{x,y}A\pare{y}\, d\nu\pare{y}}\le C\]
    uniformly, because the product of $t^{1-n}$ and $t^{n-1}$ is bounded and the uniform remainder integrates against a finite measure. It then follows that $Q_tS\ge -C\omega$, which implies 
    \[\pot P_tu=Q_tS-C_0\omega\ge -C_1\omega.\qedhere\]
\end{proof}

\begin{lem}\label{heatconverge}
    Let $K\Subset V\Subset U$, with $V$ open. Suppose $a_j\in L^1\pare{X,\omega^n}$ and $g_j,g\in C^{0,\beta}\pare{X}$ for some $0<\beta<1$ satisfy
    \[
    a_j=g_j\text{ a.e. on }V,\qquad \sup_j\norm{g_j}_{C^{0,\beta}\pare{X}}<\infty,\qquad g_j\to g\text{ uniformly on }X.
    \]
    For $t_j\searrow 0$, assume 
    \[
    t_j^{-n-1/2}e^{-c/t_j}\norm{a_j-g_j}_{L^1\pare{X}}\to 0.
    \]
    Then $P_{t_j}a_j\to g$ uniformly and $\sup_j\norm{P_{t_j}a_j}_{C^{0,\beta}\pare{K}}<\infty$.
\end{lem}
\begin{proof}
    By the triangle inequality, it suffices to show
    \[
    \norm{P_{t_j}\pare{a_j-g_j}}_{L^\infty} + \norm{P_{t_j}g_j-g_j}_{L^\infty}+\norm{g_j-g}_{L^\infty}\to 0
    \]
    and
    \[
    \norm{P_{t_j}g_j}_{C^{0,\beta}}+\norm{P_{t_j}\pare{a_j-g_j}}_{C^{0,\beta}}
    \]
    is uniformly bounded. 

    First, choose an open neighborhood $W$ with $K\Subset W\Subset V$ and consider $P_{t_j}\pare{a_j-g_j}$. Since $a_j-g_j=0$ a.e. on $V$, the kernel and its first derivative are evaluated only at pairs in $W\times\pare{X\setminus V}$, whose separation is positive. Therefore, Lemma \ref{heatbounded} gives
    \[
    \norm{P_{t_j}\pare{a_j-g_j}}_{C^1\pare{W}}\le Ct_j^{-n-1/2}e^{-c/t_j}\norm{a_j-g_j}_{L^1\pare{X}}\to 0.
    \]
    This gives convergence to zero in $C^{0,\beta}\pare{K}$. For nearby points in $K$, integrate the gradient along a short geodesic contained in $W$; for points further apart, use the supremum bound. We also have that $P_{t_j}a_j\to g$ uniformly on $K$.

    The Gaussian bound and $P_{t_j}1=1$ give 
    \[
    \norm{P_{t_j}g_j-g_j}_{L^\infty}\le C\sqbr{g_j}_\beta t^{\beta/2}.
    \]
    Moreover, differentiating $\int_X K_{t_j}\pare{x,y}\,\omega^n\pare{y}=1$ in $x$ gives
    \[
    \nabla P_{t_j} g_j\pare{x}=\int_X\nabla_x K_{t_j}\pare{x,y}\pare{g_j\pare{y}-g_j\pare{x}}\, \omega^n\pare{y}.
    \]
    Combined with the gradient bound of Lemma \ref{heatbounded}, we find
    \[
    \norm{\nabla P_{t_j}g_j}_{L^\infty}\le C\sqbr{g_j}_\beta t^{\pare{\beta-1}/2}.
    \]

    We now uniformly bound $\sqbr{P_{t_j}g_j}_\beta$. Let $x,y\in X$. If $d\pare{x,y}\le\sqrt{t_j}$, then the gradient estimate above gives
    \[
    \abs{P_{t_j}g_j\pare{x}-P_{t_j}g_j\pare{y}}\le C\sqbr{g_j}_\beta t^{\pare{\beta-1}/2}d\pare{x,y}\le C\sqbr{g_j}_\beta d\pare{x,y}^\beta.
    \]
    If instead $d\pare{x,y}\ge \sqrt{t_j}$, then comparing $P_{t_j}g_j$ to $g_j$ at both points and using $t^{\beta/2}\le d\pare{x,y}^\beta$ gives the same bound. Thus, $\sqbr{P_{t_j}g_j}_\beta\le C\sqbr{g_j}_\beta$, which is uniformly bounded. Combining everything proves the lemma.
\end{proof}


\nocite{*} 

\begin{thebibliography}{99}

\bibitem{Aronson}
Aronson, D. G.:
\emph{Bounds for the fundamental solution of a parabolic equation},
Bull. of the Amer. Math. Soc. \textbf{73} (1967), 890--896.
\href{https://doi.org/10.1090/S0002-9904-1967-11830-5}{doi:10.1090/S0002-9904-1967-11830-5}.

\bibitem{BT}
Bedford, E. and Taylor, B. A.:
\emph{A new capacity for plurisubharmonic functions},
Acta Math. \textbf{149} (1982), 1--40.
\href{https://doi.org/10.1007/BF02392348}{doi:10.1007/BF02392348}.

\bibitem{BGV}
Berline, N. and Getzler, E. and Vergne, M.:
\emph{Heat Kernels and Dirac Operators},
Grundlehren Text Editions, Springer, Berlin (2004).
\href{https://doi.org/10.1007/978-3-642-58088-8}{doi:10.1007/978-3-642-58088-8}.

\bibitem{BBGZ}
Berman, R. J. and Boucksom, S. and Guedj, V. and Zeriahi, A.:
\emph{A variational approach to complex Monge-Amp\`ere equations},
Publ. Math. de l'IH\'ES, \textbf{117} (2013), 179--245.
\href{https://doi.org/10.1007/s10240-012-0046-6}{doi:10.1007/s10240-012-0046-6}.

\bibitem{BermanDemailly}
Berman, R. J. and Demailly, J.-P.:
\emph{Regularity of plurisubharmonic upper envelopes in big cohomology classes},
in Perspectives in Analysis, Geometry, and Topology, Progr. Math. \textbf{296}, Birkh\"auser, Springer, New York, (2012) , 39--66.
\href{https://doi.org/10.1007/978-0-8176-8277-4_3}{doi:10.1007/978-0-8176-8277-4\_3}.

\bibitem{Boucksom}
Boucksom, S.:
\emph{Divisorial Zariski decompositions on compact complex manifolds},
Ann. Sci. \'Ec. Norm. Sup\'er. \textbf{37} (2004), no.~1, 45--76.
\href{https://doi.org/10.1016/j.ansens.2003.04.002}{doi:10.1016/j.ansens.2003.04.002}.

\bibitem{BEGZ}
Boucksom, S., Eyssidieux, P., Guedj, V. and Zeriahi, A.:
\emph{Monge-Amp\`ere equations in big cohomology classes},
Acta Math. \textbf{205} (2010), no.~2, 199--262.
\href{https://doi.org/10.1007/s11511-010-0054-7}{doi:10.1007/s11511-010-0054-7}.

\bibitem{CKNS}
Caffarelli, L., Kohn, J. J., Nirenberg, L. and Spruck, J.:
\emph{The Dirichlet problem for nonlinear second order elliptic equations, II: Complex Monge-Amp\`ere, and uniformly elliptic equations},
Comm. Pure Appl. Math., \textbf{38} (1985), 209--252. 
\href{https://doi.org/10.1002/cpa.3160380206}{doi: 10.1002/cpa.3160380206}.

\bibitem{DDNL}
Darvas, T. and Di Nezza, E. and Lu, C. H.:
\emph{Log-concavity of volume and complex {Monge-Amp{\`e}re} equations with prescribed singularity},
Math. Ann. \textbf{379} (2021), 95--132.
\href{https://doi.org/10.1007/s00208-019-01936-y}{doi:10.1007/s00208-019-01936-y}.

\bibitem{Demailly}
Demailly, J.-P.:
\emph{Regularization of closed positive currents and intersection theory},
J. Algebraic Geom. \textbf{1} (1992), 361--409.

\bibitem{DemaillyDinewGuedjHiepKolodziejZeriahi}
Demailly, J.-P. and Dinew, S. and Guedj, V. and Hiep, P. H. and Kołodziej, S. and  Zeriahi, A.:
\emph{Hölder continuous solutions to Monge-Ampère equations},
J. Eur. Math. Soc.  \textbf{16} (2014), no.~4, 619--647. \href{http://doi.org/10.4171/jems/442}{doi:10.4171/jems/442}.

\bibitem{DemaillyPaun}
Demailly, J.-P. and P\u{a}un, M.: 
\emph{Numerical characterization of the
K\"ahler cone of a compact K\"ahler manifold}, Ann. Math. \textbf{159} (2004), 1247--1274. \href{htpps://doi.org/10.4007/annals.2004.159.1247}{doi:10.4007/annals.2004.159.1247}


\bibitem{DNGG}
Di Nezza, E. and Guedj, V. and Guenancia, H.:
\emph{Families of singular K\"ahler--Einstein metrics},
J. Eur. Math. Soc. \textbf{25} (2023), no.~7, 2697--2762.
\href{https://doi.org/10.4171/JEMS/1249}{doi:10.4171/jems/1249}.

\bibitem{DiNezzaGuedjLu}
Di Nezza, E. and Guedj, V. and Lu, C. H.:
\emph{Geometric smoothing by the K\"ahler-Ricci flow},
(2026). \href{http://arxiv.org/abs/2607.08325}{arXiv:2607.08325}.

\bibitem{DinezzaLu}
Di Nezza, E. and Lu, C. H.:
\emph{Complex {Monge-Amp\`ere} equations on quasi-projective varieties},
J. Reine Angew. Math. \textbf{727} (2017), 145--167.
\href{http://dx.doi.org/10.1515/crelle-2014-0090}{doi:10.1515/crelle-2014-0090}.

\bibitem{DiNezzaTrapani}
Di Nezza, E. and Trapani, S.:
\emph{The regularity of envelopes},
Ann. Sci. \'Ec. Norm. Sup\'er. \textbf{57} (2024), 1347--1370. \href{https://doi.org/10.24033/asens.2592}{doi:10.24033/asens.2592}.

\bibitem{DinewGuedjZeriahi}
Dinew, S. and Guedj, V. and Zeriahi, A.:
\emph{Open problems in pluripotential theory},
Complex Var. Elliptic Equ.  \textbf{61} (2016), 902--930. \href{https://dx.doi.org/10.1080/17476933.2015.1121481}{doi:10.1080/17476933.2015.1121481}.

\bibitem{EGZ1}
Eyssidieux, P., Guedj, V. and Zeriahi, A.:
\emph{Singular K\"ahler-Einstein metrics},
J. Amer. Math. Soc. \textbf{22} (2009), no.~3, 607--639.
\href{https://doi.org/10.1090/S0894-0347-09-00629-8}{doi:10.1090/S0894-0347-09-00629-8}.

\bibitem{EGZ2}
Eyssidieux, P., Guedj, V. and Zeriahi, A.:
\emph{Viscosity solutions to degenerate complex monge-amp\`ere
equations},
Comm. Pure Appl. Math. \textbf{64} (2011), 1059--1094. \href{https://doi.org/10.1002/cpa.20364}{doi:10.1002/cpa.20364}.

\bibitem{Gaveau}
Gaveau, B.:
\emph{M\'ethodes de contr\^ole optimal en analyse complexe. I. R\'esolution d'\'equations de Monge-Amp\`ere},
J. Funct. Anal.,
\textbf{25} (1977), 391--411.
\href{https://doi.org/10.1016/0022-1236(77)90046-5}{doi:10.1016/0022-1236(77)90046-5}.

\bibitem{Grigor'yan1}
Grigor'yan, A.:
\emph{Gaussian upper bounds for the heat kernel on arbitrary manifolds},
J. Differential. Geom., \textbf{45} (1997), 33--52.
\href{https://doi.org/10.4310/jdg/1214459753}{doi:10.4310/jdg/1214459753}.

\bibitem{GuanLi}
Guan, B. and Li, Q.:
\emph{Complex Monge-Amp\`ere equations and totally real submanifolds},
Adv. Math. \textbf{225} (2010), no.~3, 1185--1223.
\href{https://doi.org/10.1016/j.aim.2010.03.019}{doi:10.1016/j.aim.2010.03.019}.

\bibitem{GuedjZeriahi1}
Guedj, V. and Zeriahi, A.:
\emph{Stability of solutions to complex Monge-Amp\`ere equations in big cohomology classes}, Math. Res. Lett., \textbf{19} (2012), 1025--1042.
\href{https://doi.org/10.4310/MRL.2012.v19.n5.a6}{doi:10.4310/MRL.2012.v19.n5.a6}.

\bibitem{HLP}
Harvey, F. R. and Lawson, H. B. Jr. and Pliś, S.:
\emph{The Richberg technique for subsolutions}, Commun. Anal. Geom., \textbf{28} (2020), 1787--1806.
\href{https://doi.org/10.4310/CAG.2020.v28.n8.a2}{doi:10.4310/CAG.2020.v28.n8.a2}.

\bibitem{Hsu}
Hsu, E. P.:
\emph{Estimates of derivatives of the heat kernel on a compact Riemannian manifold},
Proc. Amer. Math. Soc. \textbf{127} (1999), no.~12, 3739--3744.
\href{https://doi.org/10.1090/S0002-9939-99-04967-9}{doi:10.1090/S0002-9939-99-04967-9}.

\bibitem{kolodziej}
Ko{\l}odziej, S.:
\emph{The complex Monge-Amp\`ere equation},
Acta Math. \textbf{180} (1998), no.~1, 69--117.
\href{https://doi.org/10.1007/BF02392879}{doi:10.1007/BF02392879}.

\bibitem{kolodziej2}
Ko{\l}odziej, S.:
\emph{The complex Monge-Amp\`ere equation and pluripotential theory},
Mem. Amer. Math. Soc. \textbf{178} (2005), no.~840, x+64.
\href{https://doi.org/10.1090/memo/0840}{doi:10.1090/memo/0840}.

\bibitem{KolodziejNguyen}
Ko{\l}odziej, S. and Nguyen, N.-C.:
\emph{Stability and regularity of solutions of the Monge-Amp\`ere equation on Hermitian manifolds},
Adv. Math. \textbf{346} (2019), 264--304.
\href{https://doi.org/10.1016/j.aim.2019.02.004}{doi:10.1016/j.aim.2019.02.004}.

\bibitem{PS}
Phong, D. H. and Sturm, J.:
\emph{Lectures on stability and constant scalar curvature},
in Curr. Dev. Math. 2007, Int. Press, Somerville, MA (2009), 101--176.
\href{https://arxiv.org/abs/0801.4179}{arXiv:0801.4179}.

\bibitem{PS2}
Phong, D. H. and Sturm, J.:
\emph{The Dirichlet problem for degenerate complex Monge-Amp\`ere equations},
Comm. Anal. Geom. \textbf{18} (2010), no.~1, 145--170.
\href{https://doi.org/10.4310/CAG.2010.v18.n1.a6}{doi:10.4310/CAG.2010.v18.n1.a6}.

\bibitem{Sheu}
Sheu, S.-J.:
\emph{
Some Estimates of the Transition Density of a Nondegenerate Diffusion Markov Process},
Ann. Probab. \textbf{19} (1991), no.~2, 538--561.
\href{https://doi.org/10.1214/aop/1176990440}{doi:10.1214/aop/1176990440}.

\bibitem{SongTian}
Song, J. and Tian, G.:
\emph{
The K\"ahler-Ricci flow through singularities},
Invent. Math. \textbf{207} (2017), 519--595.
\href{https://doi.org/10.1007/s00222-016-0674-4}{doi:10.1007/s00222-016-0674-4}.

\bibitem{Wang}
Wang, Y.:
\emph{On the $C^{2,\alpha}$-regularity of the complex Monge-Amp\`ere equation},
Math. Res. Lett. \textbf{19} (2012), no.~4,
939--946.
\href{https://doi.org/10.4310/MRL.2012.v19.n4.a16}{doi:10.4310/MRL.2012.v19.n4.a16}.

\bibitem{Yau1978}
Yau, S.-T.:
\emph{On the Ricci curvature of a compact K\"ahler manifold and the complex Monge-Amp\`ere equation, I},
Comm. Pure Appl. Math. \textbf{31} (1978), no.~3, 339--411.
\href{https://doi.org/10.1002/cpa.3160310304}{doi:10.1002\allowbreak/cpa.3160310304}.

\end{thebibliography}
\end{document}